\documentclass[11pt,fleqn]{article}
\pdfoutput=1
\usepackage{verbatim}
\usepackage{makeidx}
\usepackage{verbatim}

\usepackage{pgf,tikz}
\usepackage{mathrsfs}
\usetikzlibrary{arrows}
\usepackage{amssymb,amsmath,graphicx}
\usepackage{epstopdf}

\newcommand{\nz}{{\mathbb{N}}}
\newcommand{\rz}{{\mathbb{R}}}
\newtheorem{axiom}{Axiom}[section]

\newtheorem{corollary}[axiom]{Corollary}
\newtheorem{definition}[axiom]{Definition}

\newtheorem{lemma}[axiom]{Lemma}
\newtheorem{observation}[axiom]{Observation}

\newtheorem{theorem}[axiom]{Theorem}
\newcommand{\proof}{\emph{Proof.}\ \ }
\newcommand{\qed}{\hfill~~$\Box$}

\newcommand{\aaa}{\alpha}
\newcommand{\bbb}{\beta}

\newcommand{\seq}[1]{\langle #1\rangle}

\begin{document} 
\sloppy

%
%
%
%
%


\title{{\bf New special cases of the Quadratic Assignment Problem with diagonally structured coefficient matrices}} 
\author{Eranda \c{C}ela
  \thanks{{\tt cela@opt.math.tu-graz.ac.at}.
  Institut f\"ur Diskrete Mathematik, TU Graz, Steyrergasse 30, A-8010 Graz, Austria}
  \and Vladimir Deineko
  \thanks{{\tt V.Deineko@warwick.ac.uk}.
  Warwick Business School, Coventry, CV4 7AL, United Kingdom}
  \and Gerhard J.\ Woeginger
  \thanks{{\tt gwoegi@win.tue.nl}.
  Department of Mathematics and Computer Science, TU Eindhoven, P.O.\ Box 513, 
  5600 MB Eindhoven, Netherlands} 
}
\date{}
\maketitle

\begin{abstract}
We consider new polynomially solvable cases of the well-known Quadratic
Assignment Problem involving coefficient matrices with  a special
diagonal structure. By combining the new special cases with polynomially
solvable special cases 
known in the literature we obtain a new and larger class of polynomially solvable special cases
of the QAP where  
one of the two coefficient  matrices involved is a Robinson matrix  with an additional structural property: 
this matrix can be represented as a conic combination of cut matrices in a certain normal form. 
The other matrix is a conic combination of a monotone anti-Monge matrix and a
down-benevolent Toeplitz matrix. We  consider the recognition
problem for the special class of Robinson matrices  mentioned above and show
that it can be solved in polynomial time. 

\end{abstract}

\section{Introduction}
\label{sec:intro}
In this paper we investigate the Quadratic Assignment Problem (QAP), which is a well-known problem
in combinatorial optimization; we refer the reader to the book \cite{Cela-book} by \c{C}ela and the 
book \cite{Burkard-book} by Burkard, Dell'Amico \& Martello for comprehensive surveys on the QAP.
The QAP in Koopmans-Beckmann form \cite{KoBe1957} takes as input two $n\times n$ square matrices
$A=(a_{ij})$ and $B=(b_{ij})$ with real entries.
The goal is to find a permutation $\pi$ that minimizes the objective function
\begin{equation}
\label{eq:qap}
Z_{\pi}(A,B) ~:=~ \sum_{i=1}^n\sum_{j=1}^n~ a_{\pi(i)\pi(j)} \, b_{ij}.
\end{equation}
Here $\pi$ ranges over the set $S_n$ of all permutations of $\{1,2,\ldots,n\}$.
In general, the QAP is extremely difficult to solve and hard to approximate.
One branch of research on the QAP concentrates on the algorithmic behavior of strongly structured
special cases; see for instance Burkard \& al \cite{BCRW1998}, Deineko \& Woeginger \cite{DeWo1998},
\c{C}ela \& al \cite{Cela2011},  \c{C}ela, Deineko \&
Woeginger~\cite{Cela2012}, or Laurent and Seminaroti~\cite{Laurent2015} 
for typical results
in this direction. 

In our paper we follow recent developments and represent several new results in
this exciting area of research.
In particular we discuss two new polynomially solvable special cases of the QAP
involving diagonally structured matrices,  the so-called
down-benevolent QAP and  the up-benevolent QAP.
 Further we  focus on the so-called combined special cases of the
  problem. They 
  arise as a  combination of different polynomially solvable special cases of
  the QAP  which involve  one coefficient
  matrix of a common type, respectively. This approach is interesting because
  it allows the identification of  new and more  complex polynomially solvable special cases
  of the QAP. 
A short discussion on one of the combined special cases presented  in  this paper was published in \cite{Cela2015a}.
\section{Preliminaries and definitions}
The are a number of results known on polynomially solvable special cases of the QAP where the coefficent matrices $A$ and $B$
possess specific structural properties. In an effort to classify the structural properties which seem to lead to polynomially 
solvable special cases of the QAP we distinguish monotonicity properties, diagonal structural properties, block structural properties 
and  properties related to so-called four-point condiditons. 
In the following we will define some matrix  classes having properties of the type mentioned above and  recall some 
 results on  
polynomially  solvable special cases of the QAP known in the literature  where these matrix classes are involved. 
 \smallskip

\begin{definition}{\bf Monotonicity properties}\\
A symmetric matrix $A=(a_{ij})$ is a \emph{Robinsonian dissimilarity} or briefly a \emph{Robinson matrix}, if for all $i<j<k$ it satisfies 
the conditions $a_{ik}\ge \max \{a_{ij},a_{jk}\}$; in words, the entries in the matrix are placed in 
non-decreasing order in each row and column when moving away from the main diagonal.

A symmetric matrix $A=(a_{ij})$ is a  \emph{Robinsonian similarity}, if for all $i<j<k$ it satisfies 
the conditions $a_{ik}\le \max \{a_{ij},a_{jk}\}$. 

An $n\times n$ matrix $B=(b_{ij})$ is called \emph{monotone}, if $b_{ij}\le b_{i,j+1}$ and $b_{ij}\le b_{i+1,j}$ 
holds for all $i,j$, that is, if the entries in every row  and column are sorted non-decreasingly from the left to the right  and from the top to the bottom, respectively.
\end{definition}

In some  QAP special cases considered in this paper  the diagonal elements of
the coefficient matrices do not impact the optimal solution. In these cases we
assume them to be zero and  set $a_{ii}=0$, for all $i$. 

The Robinson matrices were first introduced by Robinson~\cite{Robinson} in 1951
in the context of   an analysis of archaeological data. Since then they have been widely used in combinatorial data analysis; 
see the books \cite{Hubert1,Hubert2,Mirkin1,Mirkin2} and the surveys \cite{Barthelemy,Brucker} for examples 
of various applications of Robinsonian structures in quantitative psychology, analysis of DNA sequences, 
cluster analysis, etc. Special cases of the QAP involving Robinson matrices
are discussed in Laurent and Seminaroti~\cite{Laurent2015}.
\smallskip

\begin{definition}{\bf Diagonal structural properties}\\
 An $n\times n$ matrix $B=(b_{ij})$ is called a \emph{Toeplitz matrix} if it has constant entries along each of its 
diagonals;
in other words, there exists a function $f\colon \{-n+1,-n+2,\ldots,-1,0,1,\ldots,n-1\}\to \rz$
 such that $b_{ij}=f(i-j)$, for all $1\le i,j\le n$.  The  Toeplitz matric $B$ is fully determined by the function $f$ and therefore $f$ will be called {\sl the generating function of $B$}. 
If $f(i)=f(i-n)$  holds for every $i\in \{1,2,\ldots,n-1\}$,  the Toeplitz matrix $B$ is called a \emph{circulant matrix}.

A symmetric $n\times n$ Toeplitz matrix whose generating function $f$ fulfills $f(0)=0$ and $f(1)\ge f(2)\ge\ldots\ge f(n-1)$ will be called a \emph{simple
Toeplitz matrix}. (These matrices were introduced by Laurent and Seminaroti~\cite{Laurent2015}).

A symmetric $n\times n$ circulant matrix $B$  whose generating function $f$ fulfills $f(0)=0$,  
$f(1)\ge f(2)\ge \ldots\ge f(\lceil \frac{n-1}{2}\rceil)$ and $f(i) = f(n-i)$ for all $i> \lceil \frac{n-1}{2}\rceil$, is  called a \emph{DW-Toeplitz matrix}
(see Deineko and Woeginger~\cite{DeWo1998}).

A symmetric $n\times n$ Toeplitz matrix $B$  whose generating function $f$ fulfills $f(0)=0$,   
$f(1)\le f(2)\le\ldots\le f(\lceil \frac{n-1}{2}\rceil)$ and $f(i) \le f(n-i)$, 
for all $i\le \lceil \frac{n-1}{2}\rceil$, 
is called an \emph{up-benevolent  Toeplitz matrix}. 
(These matrices where introduced in~\cite{BCRW1998} as 
{\sl benevolent  Toeplitz matrices}.

Analogously a symmetric $n\times n$ Toeplitz matrix $B$ whose generating function $f$ fulfills $f(0)=0$,
 $f(1)\ge f(2)\ge\ldots\ge f(\lceil \frac{n-1}{2}\rceil)$ and $f(i) \ge f(n-i)$,  
for all $i\le \lceil \frac{n-1}{2}\rceil$,  is called a \emph{down-benevolent  Toeplitz matrix}.

Finally the attributes down-benevolent and up-benevolent will be also used for the  generating 
functions  of the Toeplitz matrices having the corresponding  properties, respectively. 
So we will talk about {\sl down-benevolent functions} and {\sl up-benevolent functions} 
defined over $\{-n+1,\ldots,-1,0,1,\ldots, n-1\}$. 
\end{definition}

The structures introduced above appeared in several special cases of the QAP
dealt with in the papers already cited in Introduction. One of the most recent 
results was presented by  Laurent \& Seminaroti in  \cite{Laurent2015}, and
will be of special interest in the context of the  paper at hand. 
In   \cite{Laurent2015} it was shown that  the $QAP(A,B)$, 
where  $A$ is a Robinson matrix and   $B$ is a  simple Toeplitz matrix is
solved to optimality by the  identity permutation.
\medskip

\begin{figure}
\unitlength=1cm
\begin{center}
\begin{picture}(10.5,5)
\put(0,0.5){\begin{picture}(3.8,4)
\includegraphics[scale=1]{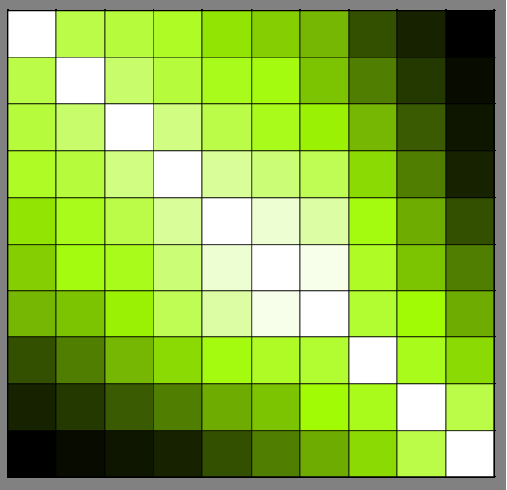}
\end{picture} }
\put(2.3,0){\makebox(0,0)[cc]{$A$}}
\put(6.5,0.5){\begin{picture}(4.2,3.6)
\includegraphics[scale=0.91]{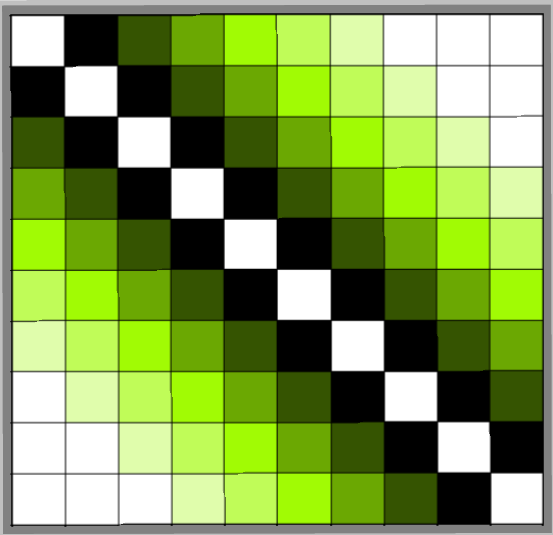}
\end{picture} }
\put(8.8,0){\makebox(0,0)[cc]{$B$}}
\end{picture}
\end{center}
\caption{Illustration to the Laurent \& Seminaroti QAP\cite{Laurent2015}: A -
  Robinsonian dissimilarity, B - simple Toeplitz matrix; the darker the colour
  the larger the entries of the matrix.}
\label{fig:RobToeplitz}
\end{figure}

To help readers to better understand structures involved in various QAP special
cases, 
we use here a color coding to visualise these
structures. Figure~\ref{fig:RobToeplitz} illustrates Robinson matrices 
 and simple Toeplitz matrices - the darker the  color the larger the  value of
 the corresponding matrix entries; 
the white colour corresponds to zero entries. 
The instances of matrices used for the illustrations can be found in Appendix.

\begin{definition}{\bf Block structural properties}\\ 
Let  a $q\times q$ matrix $P=(p_{ij})$ be fixed.  An $n\times n$ matrix $B=(b_{ij})$ is called   a
\emph{block matrix with block pattern $P$} if the following holds
\begin{itemize}
\item[(i)] there exists a partition of the set of row and column indices  $\{1,\ldots,n\}$ into $q$ 
(possibly empty) sets $I_1,\ldots,I_q$ such that for $1\le k\le q-1$ all elements of 
 $I_k$ are smaller than all
elements of     $I_{k+1}$,
\item[(ii)] for all pairs of indices $(i,j)$ with    $i\in I_k$ and $j\in I_{\ell}$ the equality  
$b_{ij}=p_{k\ell}$ holds, for all $k,\ell\in \{1,2,\ldots,q\}$
\end{itemize}
The sets $I_1,\ldots,I_q$ are called \emph{row and column blocks of matrix $B$}. 
If it is clear from context we will sometimes refer to these sets as the \emph{blocks} of matrix B.

A \emph{cut matrix} $B$ is a block matrix whose block pattern has $0$'s along the main diagonal and
$1$'s everywhere else.
A cut matrix is in \emph{CDW normal form}, if its block sizes are in non-decreasing order, i.e.
$|I_1|\le|I_2|\le\cdots\le|I_q|$ holds. (These matrices were introduced in \cite{Cela2015}.)
\end{definition}

It is easy to see that any cut matrix is a Robinson matrix. 
So it follows from \cite{Laurent2015} that if $A$ is a cut matrix, and $B$ is a simple Toeplitz matrix 
(as defined above), then the QAP is solved by the identity permutation. 
On the other hand \c{C}ela, Deineko \& Woeginger \cite{Cela2015}  have shown that if  $A$ is a cut 
matrix in CDW normal form and matrix $B$ is a monotone anti-Monge matrix (see
the definition below), then the $QAP(A,B)$ is again solved  by the identity
permutation.
\smallskip

\begin{figure}[htb]
\unitlength=1cm
\begin{center}
\begin{picture}(10.5,5)
\put(0,0.5){\begin{picture}(3.8,4)
\includegraphics[scale=0.93]{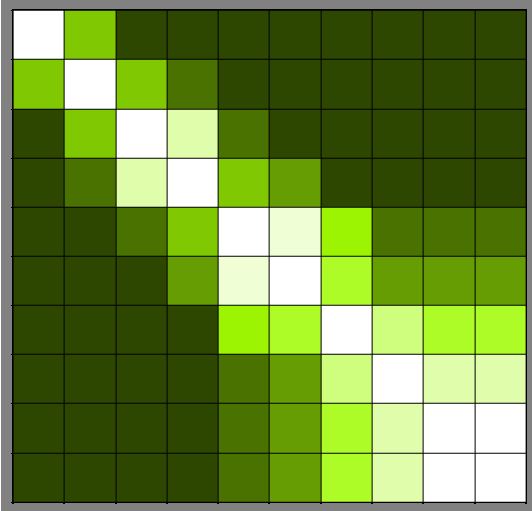}
\end{picture} }
\put(2.3,0){\makebox(0,0)[cc]{$A$}}
\put(6.5,0.5){\begin{picture}(4.2,3.6)
\includegraphics[scale=0.9]{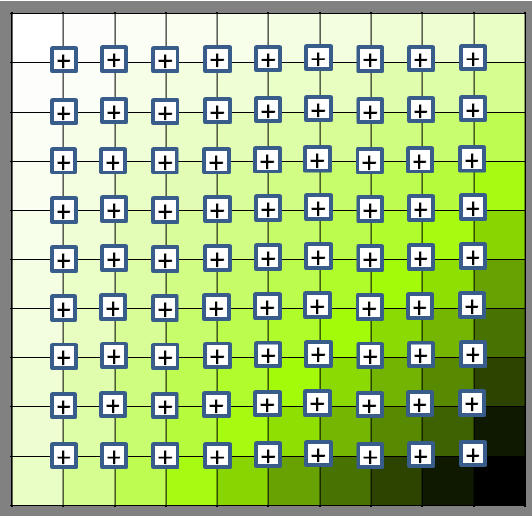}
\end{picture} }
\put(8.8,0){\makebox(0,0)[cc]{$B$}}
\end{picture}
\end{center}
\caption{Illustration to a generalisation of the Cela, Deineko \& Woeginger QAP \cite{Cela2015}: A - a conic combination of Block matrices in CDW normal form, B - Anti-Monge monotone matrix; the darker the colour
  the larger the entries of the matrix.}
\label{fig:BlockAntiM}
\end{figure}

As a  consequence  of this  result the $QAP(A,B)$  
where  $A$  is  a Robinson matrix obtained as a \emph{conic} combination of cut matrices in 
CDW normal form, 
(i.e.,  $A$ is  a linear combination of such matrices with non-negative weight coefficients) 
and $B$ is a monotone anti-Monge matrix 
 is  solved by the identity permutation. This special case is illustrated in
 Figure~\ref{fig:BlockAntiM}. 
The fulfillment of the anti-Monge inequalities is  illustrated  by  the symbol ``+".
Notice that  the block structure of  matrix $A$ is not that obvious any more in
the picture. 
\bigskip

\smallskip

In the context of the special case mentioned above  the recognition of a special class of 
Robinson matrices, namely of those Robinson matrices which can be represented as conic combinations of cut matrices in CDW normal form, becomes relevant:
\begin{quote}
Given an $n \times n$ Robinson matrix, can it be represented as a conic combination of cut matrices in CDW normal form?
\end{quote}

 The solution of this non-trivial problem is discussed in Section~\ref{recconcombCDWnormal:ssec}.  
In general {\sl  the recognition problem for  a special classes
$\cal K$ of matrices}  asks whether  a given a matrix $A$  belong to the class
$\cal K$  or not. A more general question concerns {\sl the regonition of the permuted
class $\cal K$  of matrices}: for a given matrix $A$ we ask wether
two permutations $\pi_r$  and $\pi_c$ exist, such 
that the matrix which results after permuting the rows of $A$ by $\pi_r$ and
its columns   by $\pi_c$ belongs to $K$. Yet another variant of the
recognition problem for square matrices would impose the restriction
$\pi_r=\pi_c$.
\smallskip
  
\noindent Recognition problems as formulated above can be highly non-trivial.
  There
are a number of papers  dealing with recognition problems for different
(permuted) classes of matrices, especially  also for Robinson matrices~\cite{ChFaTr,LauSemi16} in the literature. 
 \medskip

As an illustrative example for the recognition of a conic combination of cut
matrices in CDW normal form we consider the following Robinson matrix:
\[
C=\left( \begin{array}{cccccc}
 0&1&2&3&3&3\\
 1&0&2&3&3&3\\
 2&2&0&2&3&3\\
 3&3&2&0&2&2\\
 3&3&3&2&0&1\\
 3&3&3&2&1&0
\end{array} \right)
\]
which is obtained as a sum of three cut-matrices $C=C_1+C_2+C_3$; 
here matrix $C_1$ has the three blocks $\{1,2,3\}$, $\{4\}$, $\{5,6\}$, matrix $C_2$ has three blocks 
$\{1,2\}$, $\{3\}$, $\{4,5,6\}$, and matrix $C_3$ has five blocks $\{1\}$, $\{2\}$, $\{3,4\}$, $\{5\}$, and $\{6\}$.
As none of these matrices above is a cut matrix in CDW normal form, there are no reasons to assume that the QAP with 
$C$ and a monotone anti-Monge matrix $B$ is solved by the identity permutation. 
Later we will show that $C$ can indeed be represented as a conic combination of cut matrices in CDW normal form, and 
hence the corresponding QAP is solved by the identity permutation.
\medskip

\smallskip

\begin{definition}{\bf Properties related to four point conditions}\\
An $n\times n$ matrix $B$ is an \emph{Monge matrix}, if its entries are non-negative and satisfy the 
Monge  inequalities
\begin{equation}
\label{eq:Monge}
b_{ij}+b_{rs} ~\le~b_{is}+b_{rj} \mbox{\qquad for $1\le i<r\le n$ and $1\le j<s\le n$.}
\end{equation}
In other words, in every $2\times2$ submatrix the sum of the entries on the main diagonal is smaller 
than 
the sum of the entries on the other diagonal. (The Monge property essentially dates back to the work of Gaspard Monge \cite{Monge} in the 18th century.)

Analogously, an $n\times n$ matrix $B$ is an \emph{anti-Monge matrix}, if its entries are non-negative 
and satisfy the anti-Monge 
inequalities
\begin{equation}
\label{eq:anti-Monge}
b_{ij}+b_{rs} ~\ge~b_{is}+b_{rj} \mbox{\qquad for $1\le i<r\le n$ and $1\le j<s\le n$.}
\end{equation}
In other words, in every $2\times2$ submatrix the sum of the entries on the main diagonal is larger than
the sum of the entries on the other diagonal.

A symmetric $n\times n$  matrix $(c_{ij})$ is called a \emph{Kalmanson matrix},
 if it satisfies the conditions
\begin{eqnarray}
\label{eq:kal1} c_{ij}+c_{kl}\le c_{ik}+c_{jl}\\
\label{eq:kal2} c_{ik}+c_{jl}\ge c_{il}+c_{jk}
\end{eqnarray}
for all $i$,$j$,$k$ and $l$ with $1\le i<j<k<l\le n$. 
(These matrices were introduced in  1975 by Kenneth Kalmanson~\cite{K}.)
\end{definition}

Much research has been done on the effects of four point conditions in
combinatorial optimization. 
Probably the first reference to the four point conditions is due to
Supnik~\cite{Su57}, while  the term was independently introduced by Quintas \& Supnick~\cite{QuSu66} and Buneman~\cite{Bun74}. 

   Monge structures play a special  role in polynomially solvable special cases of the QAP~\cite{BCRW1998,Cela2015,Cela2011}.
We refer the reader to the survey \cite{BuKlRu1996} by Burkard, Klinz \& Rudolf for more general information on
Monge and anti-Monge structures.

Kalmanson matrices play a  role  in polynomially solvable special cases of the QAP~\cite{DeWo1998}
and also in special cases of  a number of  other   combinatorial 
optimization problems as the travelling salesman problem~\cite{K},  the prize-collecting TSP 
\cite{CFSS}, the master tour problem \cite{DRW}, the Steiner tree problem \cite{KW}, 
the three-dimensional  matching problem \cite{PSW}.

\begin{figure}
\unitlength=1cm
\begin{center}
\begin{picture}(10.5,5)
\put(0,0.5){\begin{picture}(3.8,4)
\includegraphics[scale=0.9]{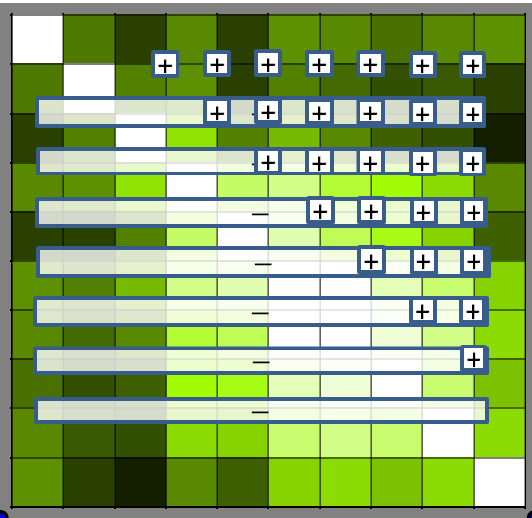}
\end{picture} }
\put(2.3,0){\makebox(0,0)[cc]{$A$}}
\put(6.5,0.5){\begin{picture}(4.2,3.6)
\includegraphics[scale=0.9]{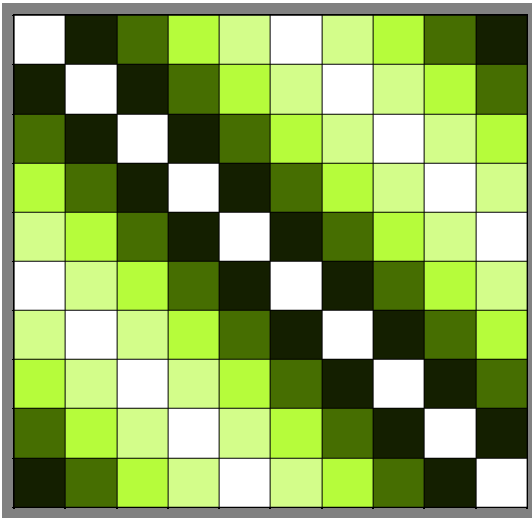}
\end{picture} }
\put(8.8,0){\makebox(0,0)[cc]{$B$}}
\end{picture}
\end{center}
\caption{Illustration to the Deineko \& Woeginger QAP\cite{DeWo1998} : A - a Kalmanson matrix, B - a DW-Toeplitz matrix.}
\label{fig:DeWo1998}
\end{figure}

Aspecial  case of the $QAP(A,B)$ involving a Kalmanson matrix  was   considered
by Deineko \& Woeginger \cite{DeWo1998}.  They showed that  the 
identity permutation is an optimal solution of the  QAP(A,B) if $A$ is a \emph{ Kalmanson} matrix and $B$ is a 
DW-Toeplitz matrix.
This special case is illustrated in Figure \ref{fig:DeWo1998}. The inequalities (\ref{eq:kal1}) and 
(\ref{eq:kal2}) fulfilled by the entries of the Kalmanson matrix $C$ are illustrated by the ``+'' 
and ``-'', respectively.

In many cases an alternative characterisation of Kalmanson matrices as
formulated in the Lemma~\ref{lemma: kalm.1} 
proved e.g.\ in \cite{DRVW,DRW} turns out to be useful.  In this
characterisation 
the fulfillment of (\ref{eq:kal2})  is  required just for quadruples of entries 
$c_{i,j}$, $c_{i+1,j+1}$, $c_{i+1,j}$ and $c_{i,j+1}$, $i,j\in  \{1,2,\ldots,n\}$, and the fulfillment of 
(\ref{eq:kal1})  is required just for quadruples of entries $c_{i,1}$, $c_{i+1,n}$, $c_{i,n}$, $c_{i+1,1}$, 
$i \in \{1,2,\ldots,n\}$.

\begin{lemma}(\cite{DRVW,DRW})
\label{lemma: kalm.1}
A symmetric $n\times n$ matrix $C$ is a Kalmanson matrix if and only if 
\begin{eqnarray}
\hspace{-1.3cm}
 c_{i,j+1}+c_{i+1,j} \le c_{ij}+c_{i+1,j+1}
&& \forall i,j: 1\le i\le n-3,~ i+2\le j\le n-1, \label{kalmx.c1}  \\
c_{i,1}+c_{i+1,n} \le c_{in}+c_{i+1,1}
&& \forall i:\ 2\le i\le n-2. \label{kalmx.c2}
\end{eqnarray}
\end{lemma}
\qed
\medskip

%

Finally we formally define (weak) sum and constant matrices.
\begin{definition}{\bf Sum matrices and constant matrices}\\
An $n\times n$ matrix $A=(a_{ij})$ is called a \emph{sum matrix}, if there exist real 
numbers $\aaa_1,\ldots,\aaa_n$ and $\bbb_1,\ldots,\bbb_n$ such that 
\begin{equation}
\label{eq:sum}
a_{ij}=\aaa_i+\bbb_j \mbox{\qquad for $1\le i,j\le n$.}
\end{equation}

An $n\times n$ matrix $A=(a_{ij})$  is called a \emph{constant matrix}, 
if all elements in the matrix are the same. 
Notice that a constant matrix is just a special case of a sum matrix.

An $n\times n$ matrix $A=(a_{ij})$ is called  a \emph{weak sum matrix}, if $A$
can be turned into a sum matrix 
by appropriately
changing the entries on its main diagonal.

An $n\times n$ matrix $A=(a_{ij})$   is a \emph{weak constant matrix}, if $A$
can be turned into a constant matrix 
by appropriately 
changing the entries on its main diagonal.
\end{definition}

We close this session with a simple but  useful observation which formalizes the
relationship between the optimal solutions of two QAP instances of the same
size, 
where the input matrices of one of them  are  obtained by permuting the
input  matrices of the other instance, respectively.  

\begin{observation}\label{permutedQAP:obse}
Let $A$ and $B$ be two $n\times n$ matrices, and 
let $\pi, \psi \in  {\cal S}_n$, where 
${\cal  S}_n$ is the set of permutations of $\{1,2,\ldots,n\}$.
Let $A^{\pi}:=(a^{\pi}_{ij})$ and $B^{\psi}:=(b^{\psi}_{ij})$ be the matrices
obtained from $A$ and $B$ by permuting them according to the permutations $\pi$
and $\psi$, respectively, i.e.\ $a^{\pi}_{ij}:=a_{\pi(i)\pi(j)}$ and
$b^{\psi}_{ij}:=b_{\psi(i)\psi(j)}$, for $1\le i,j\le n$. 
Then $Z(A^{\pi},B^{\psi},\phi)=Z(A,B,\phi\circ\pi\circ\psi^{-1})$, for all
$\phi \in {\cal S}_n$. Moreoever, if $\phi^{\ast}\in  {\cal S}_n$ is an optimal
solution of $QAP(A,B)$ then $\phi^{\ast}\circ\psi\circ\pi^{-1}$ is an optimal
solution of $QAP(A^{\pi},B^{\psi})$.  Finally, the optimal objective function
values of the two problems $QAP(A,B)$ and $QAP(A^{\pi},B^{\psi})$ coincide. 
\end{observation}
\proof
The following equalities show the first statement of the observation
\[ Z(A^{\pi},B^{\psi},\phi)=\sum_{i,j=1}^n
a^{\pi}_{\phi(i)\phi(j)}b^{\psi}_{ij}=\sum_{i,j=1}^n a_{\phi(\pi(i))
    \phi(\pi(j))}b_{\psi(i)\psi(j)}=\]
\begin{equation}
\label{obse:equ}\sum_{i^{\prime},j^{\prime}=1}^n
  a_{\phi(\pi(\psi^{-1}(i^{\prime})))\phi(\pi(\psi^{-1}(j^{\prime})))}b_{i^{\prime},j^{\prime}}=
Z(A,B,\phi\circ\pi\circ\psi^{-1})\, .
\end{equation}

The equality
above and the fact that $\phi^{\ast}$ ist an optimal solution of $QAP(A,B)$
imply that the following holds for every $\phi\in {\cal S}_n$
\[
Z(A^{\pi},B^{\psi},\phi^{\ast}\circ\psi\circ\pi^{-1})=
Z(A,B,\phi^{\ast}\circ\psi\circ\pi^{-1}\circ\pi\circ\psi^{-1})=Z(A,B,\phi^{\ast})\le \]
\[
Z(A,B,\phi)=Z(A,B,\phi\circ\psi\circ\pi^{-1}\circ\pi\circ\psi^{-1})=Z(A^{\pi},B^{\psi},\phi\circ\psi\circ\pi^{-1})
\, .
\]

Since every permutation in ${\cal S}_n$ can be written as
$\phi\circ\psi\circ\pi^{-1}$  for some $\phi \in {\cal S}_n$, the above
inequalities show that $\phi^{\ast}\circ\psi\circ\pi^{-1}$ is an optimal
solution of $QAP(A^{\pi},B^{\psi})$. 

Equality (\ref{obse:equ}) and the fact that permutation in ${\cal S}_n$ can be written as
$\phi\circ\psi\circ\pi^{-1}$  for some $\phi \in {\cal S}_n$ shows alsi that
every value of the objective function of  $QAP(A,B)$ is also a value of the
objective function of $QAP(A^{\pi},B^{\psi})$, and vice-versa. Thus the two
problems have the same set of values of the objective function and thererefore,
the same optimal value. 
\qed 
 
\smallskip
{\bf Outline of the paper.} 
In the next section we introduce new solvable cases of the QAP, the so-called
down-benevolent QAP in Section~\ref{down-benev:ssec}, and the up-benevolent QAP in
Section~\ref{up-benev:ssec}.  Then in Section~\ref{combined:ssec} we extend the variety of
known   polynomially solvable special cases of the QAP by introducing the so-called   combined
polynomially solvable  special cases.
Section~\ref{sec:cut} deals with conic representations  of  specially
 structured matrices. In Section~\ref{cutweights:ssec}   Kalmanson matrices and matrices wich are both
 Kalmanson and Robinson matrices as characterised in terms of  conic
 combinations  of particular cut
 matrices. These results are then used in Section~\ref{recconcombCDWnormal:ssec} to  give a
 characterisation of  conic
 combinations of cut-matrices in CDW normal form. This characterisation allows 
 the efficient recognition of  conic combinations of cut matrices in CDW
 normal form  which is a relevant issue because  these conic combinations are involved in a combined special
 case of the QAP (the first combined special case described
 in Section~\ref{combined:ssec}). 
The conclusions and some issues for further research conclude the paper in Section~\ref{conclu:sec}.

\section{New special cases of the  QAP solved by the identity permutation}
\subsection{The down-benevolent  QAP}\label{down-benev:ssec}

\begin{figure}
\unitlength=1cm
\begin{center}
\begin{picture}(10.5,5)
\put(0,0.5){\begin{picture}(3.8,4)
\includegraphics[scale=0.9]{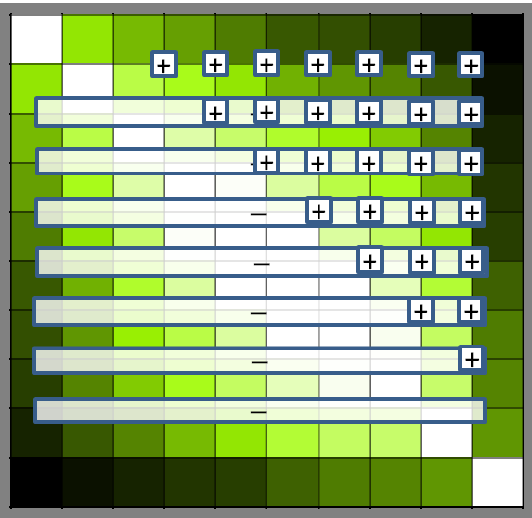}
\end{picture} }
\put(2.3,0){\makebox(0,0)[cc]{$A$}}
\put(6.5,0.5){\begin{picture}(4.2,3.6)
\includegraphics[scale=0.9]{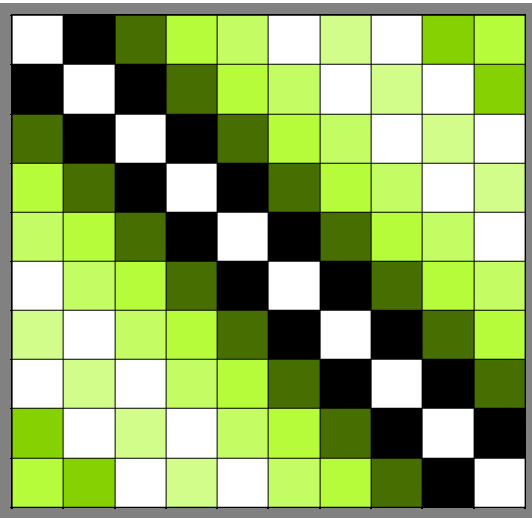}
\end{picture} }
\put(8.8,0){\makebox(0,0)[cc]{$B$}}
\end{picture}
\end{center}
\caption{Illustration to Theorem \ref{theo:Malevolent}: A - a Kalmanson and
  Robinsonian matrix, B -
 a down-benevolent Toeplitz matrix.}
\label{fig:Malevolent}
\end{figure}

In this section we consider a special case $QAP(A,B)$ which we call {\sl down-benevolent QAP}:
  $A$ is  both a Robinson matrix  and a Kalmanson matrix, and $B$ is a down-benevolent 
Toeplitz matrix. We show that this special case, illustrated in Figure~\ref{fig:Malevolent}, 
is solved by the identity permutation. 

Notice  that a  simple Toeplitz matrix is a special case of a down-benevolent Toplitz matrix. 
Analogously a DW-Toeplitz matrix is also a special case of a  down-benevolent Toplitz matrix. 
Thus, the QAP special case considered here is related to  the  QAP special cases 
considered in  \cite{Laurent2015} and in \cite{DeWo1998}. 
In \cite{Laurent2015} it was shown that the $QAP(A,B)$ with $A$ being a Robinson matrix and $B$ 
being a simple Toeplitz matrix is solved by the identity permutation.  
In \cite{DeWo1998} it was shown that the $QAP(A,B)$ with $A$ being a Kalmanson  matrix and 
$B$ being a DW-Toeplitz matrix is solved by the identity permutation. 
The new special case involves a matrix $B$ from a class which is striclty larger than  both classes 
of matrices $B$ considered in \cite{DeWo1998,Laurent2015}. 
However the matrix $A$  is required to have more restrictive properties than in
\cite{DeWo1998,Laurent2015}: $A$ is
  both a Robinson and a Kalmanson matrix.

In the following we will work  with  some particular $0$-$1$ functions. 
 \begin{definition}\label{spec:fnc}
For  $n\in \nz$ and  $i \in \nz$, $\lceil \frac{n-1}{2} \rceil< i \le n-1$,
define a function $g_n^{(i)} \colon \{-n+1,-n+2,\ldots,-1,0,1,\ldots,n-1\}\to \{0,1\}$ with 
\[ g_n^{(i)}(x)=1 \mbox{ if } x\in \{i,n-i\} \mbox{ and }   g_n^{(i)}(x)=0 \mbox{ if } i\not\in\{i,n-i\}\, .\]
\end{definition}
It can be easily seen that every $n\times n$ down-benevolent Toeplitz matrix can be obtained 
from a DW-Toeplitz matrix by subtracting from it a conic combination of Toeplitz matrices generated 
by function $g_n^{(i)}$. More precisely the following lemma holds  
\begin{lemma}
Let $B$ be an  $n\times n$ down-benevolent Toeplit matrix. Then there exists an $n\times n$  
DW-Toeplitz matrix $B^{\prime}$ and the nonegative numbers $\beta_i$, $\lceil \frac{n-1}{2} \rceil< i \le n-1$, such that $B=B^{\prime}-\sum_{i= \lceil \frac{n-1}{2} \rceil+1}^{n-1}\beta_iT^{(i)}$,  where $T^{(i)}$ is an $0$-$1$ $n\times n$ Toeplitz matrix generated by the function $g_n^{(i)}$ for $\lceil \frac{n-1}{2} \rceil< i \le n-1$.  
\end{lemma}  
\proof
Let $f\colon \{-n+1,-1,0,1,n-1\}\to \rz$ be the generating function of $B$. 
Define a symmteric function  
$f^{\prime}\{-n+1,-1,0,1,n-1\}  \to \rz$  by $f^{\prime}(i):=f(i)$ if $|i| \le \lceil \frac{n-1}{2} \rceil$ and $f^{\prime}(-i):=f^{\prime}(i):=f^{\prime}(n-i)=f(n-i)$ if 
$i > \lceil \frac{n-1}{2} \rceil$. Let $B^{\prime}$ be a Toeplitz matrix with generating function 
$f^{\prime}$. By definition $B^{\prime}$ is  a DW-Toeplitz matrix. Define $\beta_i:=f(n-i)-f(i)\ge 0$, 
for $i\in \nz$, $ \lceil \frac{n-1}{2} \rceil< i \le n-1$. Thus $f(i)=f^{\prime}(i)-\beta_i$ holds for 
evry  $i\in \nz$, $ \lceil \frac{n-1}{2} \rceil< i \le n-1$, and this implies $B=B^{\prime}-\sum_{i= \lceil \frac{n-1}{2} \rceil+1}^{n-1}\beta_iT^{(i)}$. 
\qed
\smallskip

Let $n$ be an  arbitrary but fixed natural number  and  $i\in \nz$, $\lceil \frac{n-1}{2} \rceil <i\le n-1$.
Consider the  maximization version of the $QAP(A,T^{(i)})$ with an $n\times n$   Kalmanson
 matrix $A$ which is also a Robinson matrix,  and a $T^{(i)}$ a Toeplitz matrix
 generated by $g_n^{(i)}$  with $\lceil \frac{n-1}{2} \rceil <i\le n-1 $,
 i.e.\ the  optimization problem 
$\max \{ Z_\pi(A,T^{(i)})\colon \pi \in {\cal S}_n\}$, where 
${\cal  S}_n$ is the set of permutations of $\{1,2,\ldots,n\}$. 
  Observe that $T^{(i)}$ contains  exactly $2(n-i)$ ones placed in pairwise
  symmetric positions with respect to the diagonal.  The $1$-entries above the
  diagonal lie in the rows with indices $\{1,2,\ldots,n-i\}$ and in the columns
  with indices  $\{i+1,i+2,\ldots n\}$ with exactly one $1$-entry per row and column.
Notice that since $i>  \lceil \frac{n-1}{2} \rceil$ the sets of row indices
and column  indices above do
 not intersect. The objective function value of the $QAP(A,T^{(i)})$
 corresponding to permutation
 $\pi \in {\cal S}_n$  is given as
\[ Z_\pi(A,T^{(i)})  = \sum_{k=1}^n\sum_{j=1}^n a_{\pi(k)\pi(j)}T^{(i)}_{kj}= 
\sum_{k=1}^{n-i}a_{\pi(k)\pi(k+i)} +  \sum_{k=i+1}^{n}a_{\pi(k)\pi(k-i)}=2\sum_{k=1}^{n-i}a_{\pi(k)\pi(k+i)}\, , \]
where the last equality holds because $A$ is by definition a symmetric matrix. 
Thus  $Z_\pi(A,T^{(i)})$ is just the sum of $2(n-i)$ pairwise symmetric
(non-diagonal) entries selected from  $A$, such that in every row and colum
there is at most one selected entry. 
Notice that if  each  pair of symmetric entries is represented by the
above-diagonal entry than the goal function the $QAP(A,T^{(i)})$ can be seen as
twice the sum of $n-i$ above-diagonal entries selected in $A$ such that the row
indices of selected entriws build a set $R$, the column indices of selected
entries  build a set $C$, and  $R\cap
C=\emptyset$ as well as $|R|=|C|=n-i$ hold.

Vice-versa, consider 
a set of row indices $R$ and a set of column indices $C$ with $R\cap
C=\emptyset$, $|R|=|C|=n-i$ and a bijection $\phi\colon R\to C$. 
Now select in $A$ the entries $a_{i\phi(i)}$, for $i\in R$, together with their
symmetric counterparts. 
It can be easily seen that the overall sum of these selected  entries equals  $Z_\pi(A,T^{(i)})$ 
for  any $\pi \in {\cal S}_n$  with  $\{\pi(1),\pi(2),\pi(n-i)\}=R$, $\pi(i+j)=\phi(\pi(j))$, for $1\le j\le n-i$. 
Thus the maximization version of the $QAP(A,T^{(i)})$ of size $n$ with $i$ such that 
$\lceil \frac{n-1}{2} \rceil <i\le n-1 $ is equivalent to the following selection problem 

\begin{quote}
{\bf Selection problem}

Input: $n\in \nz$, a Kalmanson and Robinson $n\times n$ matrix $A$,  $i\in \nz$ such that $\lceil \frac{n-1}{2} \rceil <i\le n-1 $ holds.

Output: Select $(n-i)$ above-diagonal  entries $a_{r_jc_j)}$, $1\le j\le
n-i$, from  $A$, such that the averall  sum $\sum_{j=1}^i a_{r_i,c_i}$ of the
selected entries is maximized,  under the condition that  the  set $R=\{r_j\colon
1\le j\le n-i\}$  of   row indices of
the selected entries  and the set $C=\{c_j\colon 1\le j\le n-i\}$ of column
indices of the selected entries fulfill $R\cap C=\emptyset$, $|R|=|C|=n-i$.   
\end{quote}

Since  the selection problem we select $n-i$ entries at most one entry  each
ro, its  solution can be  represented by a pair $(R,\phi)$,  where $R$ is the
set of  indices 
of the selected rows and   $\phi$ is 
 injective mapping $\phi\colon R\to \{1,2\ldots,n\}$ which maps each $r\in R$
 to the column index of the entry $a_{r,\phi(r)}$ selected in row $r$. 
Then, clearly,  $R\cap C=\emptyset$ holds   with $C=\{\phi(r)\colon r\in R\}$ .  
 If an entry $a_{jl}$ is selected in a solution $(R,\phi)$, i.e.\ $phi(j)=l$,
 $j\in R$,  we will say that
row index $j$ is matched with column index $l$ and column index $l$ is matched
with row index $j$ in that solution. 

Next we show that  the maximization version of
$QAP(A,T^{(i)})$, with $n\in \nz$ and  $i\in \nz$ with $\lceil \frac{n-1}{2}
\rceil <i\le n-1$, is soved by the identity permutation.

\begin{lemma}\label{KalmRobToeplOneStripe}
The maximization version of the $QAP(A,T^{(i)})$ with an $n\times n$   Kalmanson and Robinson 
 matrix $A$ and $T^{(i)}$ a Toeplitz matrix generated by $g_n^{(i)}$ with 
$i >  \lceil \frac{n-1}{2} \rceil$ is solved to optimality by the identity permutation. 
\end{lemma}
\proof
We consider the corresponding selection problem and show that it is solved to optimality by
 selecting the entries $a_{1,1+i}$, $a_{2,2+i}$, ..., $a_{n-i,n}$. Clearly, 
 this selection is feasible and  corresponds to the identity permutation in ${\cal S}_n$ as an
 optimal solution of the maximization version of the $QAP(A,T^{(i)})$ and this woud complete the proof.

Consider an optimal solution $(R,\phi)$  of the selection problem where the   row indices
of the selected entries build the set $R=\{r_1,r_2,\ldots,r_{n-i}\}$ and 
and  $R$ the  corresponding column indices are $\phi(r_j)$, for $1\le
j\le n-i$.  Then, clearly,   $R\cap \{\phi(r_j)\colon 1\le i\le
n-i\}=\emptyset$ holds. 

 Assume w.l.o.g. that $r_1<r_2<\ldots r_{n-i}$. 
First we claim  that there exists an optimal solution with
$\max \{ r_j\colon 1\le j\le n-i\}<\min\{\phi(j)\colon 1\le j\le n-i\}$,
i.e.\ an optimal solution with the following property 

(P): any  row indices of a selected entries is smaller that any column index of a selected entry
 
Assume the optimal solution  $(R,\phi)$  above does not have  Property P
Then there exist two indices $j,l\in \{1,2,\ldots,n-i\}$ such that $\phi(r_l) <
r_j$  holds. Let $r_j$ be the smallest element in $R$ for which such a 
column index of a selected entry smaller than $r_l$ exists, i.e.\ $\phi(R)\cap
\{1,2,\ldots,r_{j-1}\}\neq\emptyset$,  and let $r_l$  be
such that  $\phi(r_l)$ is the smallest column index of a selected entry which
is smaller than $r_j$, i.e. $\phi(r_l)=\min \phi(R)\cap \{1,2,\ldots,r_{j-1}\}$. 

The we clearly have $r_l<\phi(r_l)<r_j<\phi(r_j)$. 
Consider a pair  $(R',\phi')$ obtained by {\sl exchanging} $r_j$ and $\phi(r_l)$
 in the following sense:
 \[ R':=(R\setminus \{r_j\})\cup \{\phi(r_l)\}\]
\[ \phi'(r)=\phi(r),  \mbox{ $\forall r\in R\setminus \{r_j,r_l\}$  and }
\phi'(r_l)=r_j \, \phi'(\phi(r_l))=\phi(r_j).\]

 $(R',\phi')$ is a feasible soltion of the selection problem because the
two entries $a_{r_l\phi(r_l)}$, $a_{r_j\phi(r_j)}$ selected with $(R,\phi)$
are replaces by the entries  $a_{r_lr_j}$, $a_{\phi(r_l)\phi(r_j)}$ selected with $(R',\phi')$  and the
sets $R'$, $C'$ of the row and column  indices of selected entries,
respectively, fulfill clearly the properties $R'\cap C'=\emptyset$,
$|R'|=|C'|=n-i$. Moreover  inequality \ref{eq:kal1}
in the definition of Kalmanson matrices which  applies because  $A$ is a
Kalmanson 
matrix we get:
\[ a_{\phi(r_l)r_j}+a_{\phi(r_l)\phi(r_j)}\ge
a_{r_l\phi(r_l)}+a_{r_j,\phi(r_j)}\, ,\]
and thus the solution $(R',\phi')$ is not worse than the optimal solution
$(R,\phi)$, hence it is also an optimal solution. If $(R',\phi')$ does not
have  property P, then there will be again a smallest row index $r_k$
of a selected entry for which there exists a column  index of a selected entry
which smaller than $r_k$. Notice that in this case $r_k$ has to be larger than
$r_j$ because for  indices    in the set  $(R\cup)C \cap
\{1,2,\ldots,r_{j-1}\}$ the following statement holds:  any  row index os an
entry selected by the solution $(R',\phi')$ is   smaller than
any  column index of an entry selected by  $(R',\phi')$. 
So, if $(R',\phi')$ does not have  property P, then we could perform again an
exchange to abtain a new optimal solution as described above and repeat this
step as long as the current optimal solution does not have property P. 
The proces would terminate because the smallest row index of a selected entry
for which there is an even smaller column index of a selected entry, increases
in every repetion of the exchange step described abive. 
So the claim about the existence of the optimal solution with the property P is
proven.  

Let $(R,\phi)$ be an optimal solution with property P and let
$R=\{r_1,r_2,\ldots,r_{n-i}\}$ be the row indices of the selected entries with
$r_1<r_2<\ldots r_{n-i}$.
We can assume w.l.o.g.\ that $r_l<r_j$ implies $\phi(r_l)<\phi(r_j)$, for all
$l,j \in \{1,2,\ldots,n-i\}$. Indeed if there exists a pair $r_l<r_j$ for which
$\phi(r_l)>\phi(r_j)$, then consider the solution $(R,\phi')$ with
$\phi(r_k)=\phi(r_k)$ for all $k\in \{1,2,\ldots,n-i\}\setminus \{j,l\}$ and
$\phi'(j)=\phi(l)$, $\phi'(j)=\phi(l)$. Thus the entries 
 $a_{r_l\phi(r_l)}$, $a_{r_j\phi(r_j)}$ selected with $(R,\phi)$
are replaces by the entries  $a_{r_l\phi(r_j)}$, $a_{\phi(r_j)\phi(r_l)}$
selected with $(R,\phi')$. 
From   inequality \ref{eq:kal2}
in the definition of Kalmanson matrices which  applies because  $A$ is a
Kalmanson matrix we get 
\[ a_{r_l\phi(r_j)}+a_{r_j\phi(r_l)}\ge
a_{r_l\phi(r_l)}+a_{r_j,\phi(r_j)}\, ,\]
and thus the solution $(R,\phi')$ is not worse than the optimal solution
$(R,\phi)$, and hence it is also an optimal solution.

Let us denote the set of column indices of the  entries  selected with $(R,\phi)$ by
$C=\{c_1,c_2,\ldots,c_{n-i}\}$ where $c_1<c_2<\ldots<c_{n-i}$. Then the selected
entries are  $a_{r_j,c_j}$, $j=1,2,\ldots,n-i$, and $r_1 < r_2 < \ldots <
r_{n-i} < c_1<c_2 < \ldots < c_{n-i}$ holds. Then clearly $j\le r_j$ and
$c_j\le i+j$ hold, for all $j=1,2,\ldots, n-1$. 
Since matric $A$ is a Robinson matrix the above inequalities imply
$a_{r_jc_j}\le a_{j,i+j}$, for all $j=1,2,\ldots, n-1$, 
and thus
\[   \sum_{j=1}^{n-1} a_{r_jc_j} \le \sum_{j=1}^{n-1} a_{j,i+j} \, .\]
Hence selecting the  entries $a_{1,1+i}$, $a_{2,2+i}$, ..., $a_{n-i,n}$ 
is not worse then the optimal solution $(R,\phi)$, which means that
$a_{1,1+i}$, $a_{2,2+i}$, ..., $a_{n-i,n}$ is an optimal selection and
completes the proof. 
\qed

\begin{theorem}
\label{theo:Malevolent}
The $QAP(A,B)$ where $A$ is both a Robinson matrix  and a Kalmanson matrix, and
$B$ is a down-benevolent Toeplitz matrix, 
is solved to optimality by the identity permutation.
\end{theorem}
\proof
Consider an abritrary  $n\times n$ down-benevolent Toeplitz matrix generated by a function 
$f\colon \{-n+1,-n+2,\ldots,-1,0,1,\ldots,n-1\}\to \rz$
 such that $b_{ij}=f(i-j)$, for all $1\le i,j\le n$.  By definition we have  $f(0)=0$,
 $f(1)\ge f(2)\ge\ldots\ge f(\lceil \frac{n-1}{2}\rceil)$ and $f(i) \le f(n-i)$,  
for all $i>\lceil \frac{n-1}{2}\rceil$. $B$ can be represented as
\begin{equation}\label{decomp:down_benev}
 B=B_1-\sum_{i=\lceil \frac{n-1}{2}\rceil+1}^{n-1} (f(n-i)-f(i))T^{(i)}\, , 
\end{equation}
 where
$T^{(i)}$ is a Topelitz matrix generated by the function $g_n^{(i)}$ defined in
Definition~\ref{spec:fnc}, and $B_1$ is a DW-Toeplitz matrix generated by the
function $f'  \colon \{-n+1,-n+2,\ldots,-1,0,1,\ldots,n-1\}\to \rz$ with
$f'(i)=f(i)$ for $1\le i \le \lceil \frac{n-1}{2}\rceil$ and $f'(i)=f(n-i)$ for
$i>\lceil \frac{n-1}{2}\rceil$.
\smallskip

Equation (\ref{decomp:down_benev}) implies:
\[ Z(A,B,\pi)=Z(A,B_1,\pi)-\sum_{i=\lceil \frac{n-1}{2}\rceil+1}^{n-1}
(f(n-i)-f(i))Z(A,T^{(i)})\, .\]
It was proven in \cite{DeWo1998} that $QAP(A,B_1)$ with a Kalmanson matrix $A$
and a DW-Toeplitz matrix $B_1$ is solved to optimality by the identity
permutation. Lemma~\ref{KalmRobToeplOneStripe} implies  that the
maximization version of the $QAP(A,T^{(i)})$ with $T^{(i)}$ as above  is solved to optimality by the
identity permutation for all $i> \lceil \frac{n-1}{2}\rceil$.
Summarizing we get:
\[ Z(A,B,\pi)=   Z(A,B_1,\pi)-\sum_{i=\lceil \frac{n-1}{2}\rceil+1}^{n-1}
(f(n-i)-f(i))Z(A,B,T^{(i)})\ge \]
\[Z(A,B_1,id) - \sum_{i=\lceil \frac{n-1}{2}\rceil+1}^{n-1}
(f(n-i)-f(i))Z(A,T^{(i)},id)=Z(A,B,id), \mbox{ for all $\pi \in {\cal S}_n$,}\]
where $id$ denotes the identity permutation in ${\cal S}_n$. 
Thus $id$ minimizes $Z(A,B,\pi)$ for $\pi \in {\cal S}_n$ and this completes
the proof. 
\qed

\subsection{The up-benevolent QAP}\label{up-benev:ssec}

\begin{figure}
\unitlength=1cm
\begin{center}
\begin{picture}(10.5,5)
\put(0,0.5){\begin{picture}(3.8,4)
\includegraphics[scale=0.9]{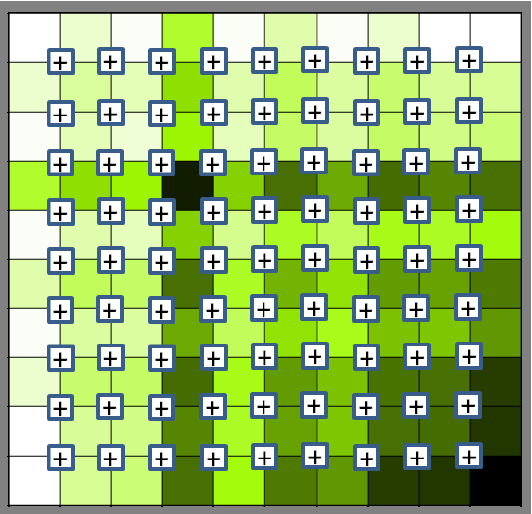}
\end{picture} }
\put(2.3,0){\makebox(0,0)[cc]{$A$}}
\put(6.5,0.5){\begin{picture}(4.2,3.6)
\includegraphics[scale=0.9]{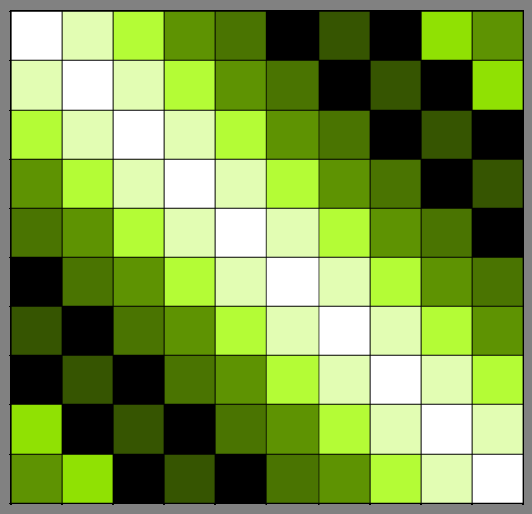}
\end{picture} }
\put(8.8,0){\makebox(0,0)[cc]{$B$}}
\end{picture}
\end{center}
\caption{Illustration to the special case of  Burkard \& al. \cite{BCRW1998} QAP: $A$ - an
  anti-Monge matrix, to be permuted with $\pi^*$; $B$ - an up-benevolent Toeplitz matrix.}
\label{fig:Bene}
\end{figure}

Burkard \& al.~\cite{BCRW1998} have considered the $QAP(A,B)$ with a monotone
anti-Monge matrix $A$ with nonnegative entries and an up-benevolent Toeplitz  $B$ (called  \emph{
  benevolent} Toeplitz matrix in the original paper). It was proven in
\cite{BCRW1998}  that the so-called Supnick permutation  $\pi^*=\seq{1,3,5,\ldots, 6,4,2}$
is an optimal solution to that QAP. An illustration of this special  case is presented in  Figure~\ref{fig:Bene}. 

In \cite{BCRW1998} it was shown that for any $n\in \nz$ the $n\times n$    monotone anti-Monge matrices with
non-negative entries form a cone whose extremal rays are given by the  $0$-$1$
matrices $R^{(p,q)}=…\big (r_{ij}^{(p,q)}\big )$, $1\le p,q\le n$, defined by
$r^{(p,q)}_{ij}=1$ for $n-p+1\le i\le n $ and $n-q+1\le j\le n$, and
$r^{(p,q)}_{ij}=1$, 
otherwise. 
 As a consequence of this fact  it can be shown that the extremal rays of the cone of {\sl  symmetric} monotone
anti-Monge  matrices with nonnegative entries are given as described by the
following lemma (see Rudolf and Woeginger~\cite{RuWo95}, Burkard et al.~\cite{BCRW1998}, and
\c{C}ela et al.~\cite{Cela2015}).
\begin{lemma}\label{SymMonAMonge:cone}
The symmetric monotone anti-Monge  matrices with nonnegative entries form a cone with extremal rays
given as $R^{(p,q)}+R^{(q,p)}$, for $1\le p <q\le n$,  and $R^{(p,p)}$,
  for $1\le p\le n$. 
\end{lemma}
Let us denote $\bar{R}^{(p,q)}:= R^{(p,q)}+R^{(q,p)}=\big
(\bar{r}^{(p,q)}_{ij}\big )$, for
$1\le p< q\le n$. The matrices $\bar{R}^{(p,q)}$ are
explicitly  given as  follows. 
$$
\bar{r}^{(p,q)}_{ij}=\left\{ \begin{array}{ll}
2, & n-p+1 \le i,j \le n\\ 
1, & n-q+1\le i \le n-p, n-p+1\le j \le n\\
1, & n-p+1 \le i \le n,\ n-q+1 \le j \le n-p \\
0, & \textrm{otherwise}
\end{array}\right.
$$
Further let us denote $\bar{R}^{(p,p)}=R^{(p,p)}$, $1\le p\le n$, for the sake
of completeness. 
\medskip

\begin{figure}
\unitlength=1cm
\begin{center}
\begin{picture}(10.5,5)
\put(0,0.5){\begin{picture}(3.8,4)
\includegraphics[scale=0.9]{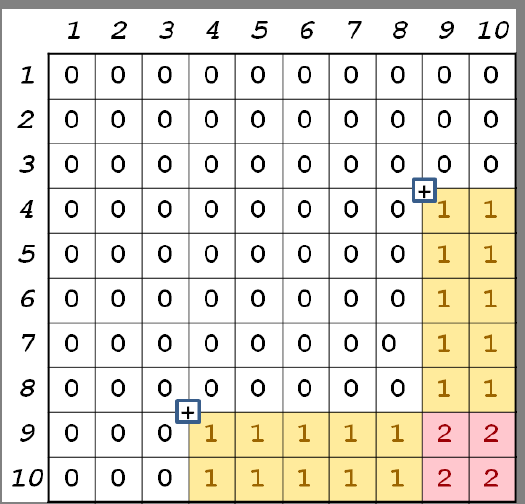}
\end{picture} }
\put(2.3,0){\makebox(0,0)[cc]{$A$}}
\put(6.5,0.5){\begin{picture}(4.2,3.6)
\includegraphics[scale=0.9]{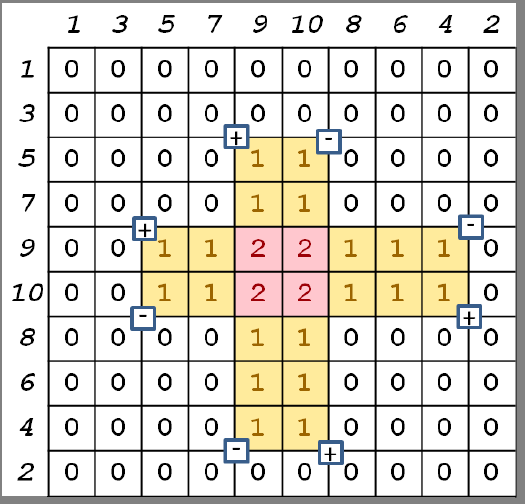}
\end{picture} }
\put(8.8,0){\makebox(0,0)[cc]{$A^{\pi^*}$}}
\end{picture}
\end{center}
\caption{Illustration to permuted matrices: $A:=\bar{R}^{(2,7)}$ - a $10\times
  10$ symmetric monotone  anti-Monge matrix; 
$A^{\pi^*}$ - a permuted anti-Monge matrix: anti-Monge inequalities are
violated.}
\label{fig:ShiftedPermuted}
\end{figure}

According to Observation~\ref{permutedQAP:obse}, if $\pi^{\ast}$ is an optimal
solution of $QAP(A,B)$ with a symmetric monotone anti-Monge matrix $A$ and  an
up-benevolent matrix $B$, then $id$ is an optimal solution of $QAP(A^{\pi^{\ast}},B)$.  
In particular this clearly holds for $A=\bar{R}^{(p,q)}$ with $1\le p\le q\le n$.  
Notice that, in general, the permuted matrix $A^{\pi^*}$ is not an anti-Monge
matrix any more. 
Figure \ref{fig:ShiftedPermuted} illustrates the effect of permuting the
$10\times 10$ matrix  $\bar{R}^{(2,7)}$ according to permutation $\pi^\ast$.
\smallskip

\smallskip

By taking  a closer look at the matrices $\bar{R}^{(p,q,\pi^{\ast})}$
obtained by permuting $\bar{R}^{(p,q)}$ according to the Supnick permutation
$\pi^{\ast}$, $1\le p\le q\le n$, we can observe that they are given as follows
\[\bar{r}^{(p,q,\pi^{\ast})}_{ij}=\left \{ \begin{array}{ll}
2 & \lceil \frac{n-p}{2}\rceil +1 \le i,j \le n-\lfloor
\frac{n-p}{2}\rfloor \\ [0.2cm]
1 & \lceil \frac{n-p}{2}\rceil - \lfloor \frac{q-p}{2}\rfloor +1 \le i\le
\lceil \frac{n-p}{2}\rceil \, ,\frac{n-p}{2}+1 \le j \le  n-\lfloor
\frac{n-p}{2}\rfloor\\[0.2cm]
1 & \lceil \frac{n-p}{2}\rceil +1 \le i\le
n-\lfloor \frac{n-p}{2}\rfloor \, ,n- \lfloor \frac{n-p}{2}\rfloor +1 \le j \le  n-\lfloor
\frac{n-p}{2}\rfloor+\lceil \frac{q-p}{2} \rceil \\[0.2cm]
0 &\mbox{otherwise}
\end{array}\right .\]
for $1\le i\le j\le n$. 
\smallskip

Thus the non-zero entries of these permuted matrices build a kind of a cross
whith entries equal to $2$ at the center of the cross and entries equal to $1$ at the arms
of the cross (see also Figure~\ref{fig:ShiftedPermuted}). Now consider a
transformation of the   matrix
$\bar{R}^{(p,q,\pi^{\ast})}$ realized by sliding  the cross of non-zero entries  along the
diagonal such that its arms do not wrap around  the border of the matrix  (they
may touch the border but should not wrap around it). This transformation is
relised by permuting (the rows and columns of) $R^{(p,q,\pi^{\ast})}$ 
according to a shift $\sigma_u\in {\cal S}_n$ of the form $\langle u,u+1,\ldots,n,1,\ldots,u-1\rangle$
 with $1< u\le \lceil \frac{n-p}{2}\rceil - \lfloor \frac{q-p}{2}\rfloor
 +1$, or 
 $n-\lfloor\frac{n-p}{2}\rfloor+\lceil \frac{q-p}{2}\rceil +1
 \le u\le n$.
\smallskip

 Let us denote by
$C^{(p,q,u)}$ the matrix obtained from $R^{(p,q,\pi^{\ast})}$ by permuting it
according to $\sigma_u$ with $u$ as described above. 
Obviously  $Z(C^{(p,q,u)},B,id)=Z(\bar{R}^{(p,q,\pi^{\ast})},B,id)$
holds for all $1\le p\le q\le n$, for all possible values  of  $u$ as given above,
and for any Toeplitz matrix $B$. This is due to the facts that a) the permutation
$\sigma_u$ shifts non-zero entries of $\bar{R}^{(p,q)}$  along lines parallel
to the main diagonal and  b) a Toeplitz matrix has constant entries along any
line parallel to the main diagonal. Combined with the third statement of
Observation~\ref{permutedQAP:obse} the above equation shows $id$ is also the
optimal solution of $QAP(C^{(p,q,u)},B)$.
This observation motivates the following definition.
\smallskip

\begin{figure}[htb]
\unitlength=1cm
\begin{center}
\begin{picture}(10.5,5)
\put(0,0.5){\begin{picture}(3.8,4)
\includegraphics[scale=0.9]{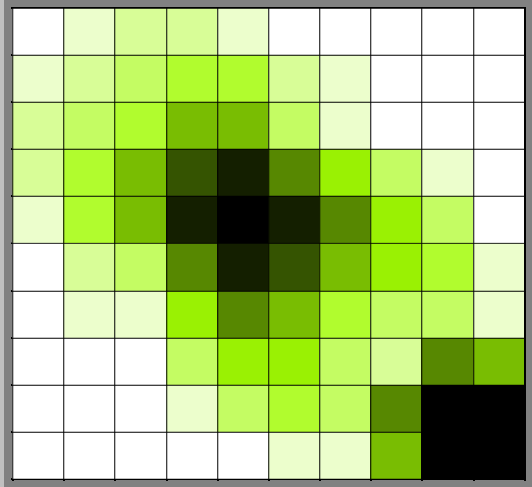}
\end{picture} }
\put(2.3,0){\makebox(0,0)[cc]{$A$}}
\put(6.5,0.5){\begin{picture}(4.2,3.6)
\includegraphics[scale=0.9]{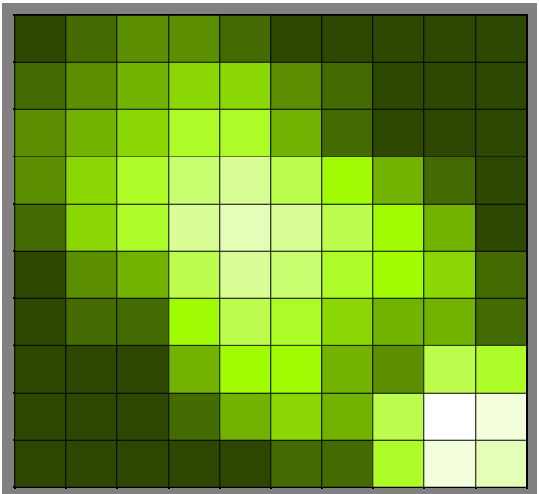}
\end{picture} }
\put(8.8,0){\makebox(0,0)[cc]{$A^-$}}
\end{picture}
\end{center}
\caption{Illustration to PS matrices: $A$, a PS monotone anti-Monge matrix;
  $A^-$, a  PS monotone Monge matrix; Monge (anti-Monge) conditions are not satisfied any more.}
\label{fig:8Rays}
\end{figure}

\begin{definition}
A symmetric $n\times n$ matrix $A^{\prime}$ is called a {\sl permuted-shifted
monotone  anti-Monge
matrix} ({\sl PS monotone anti-Monge matrix }), if it can be obtained as a conic
combination of matrices $C^{(p,q,u)}$ obtained 
$\bar{R}^{(p,q,\pi^{\ast})}$ by first permuting  acoording to $\pi^*$ and then
by applying a shift $\sigma_u$  to it, for $1\le p\le q\le n$,  $1< u\le \lceil \frac{n-p}{2}\rceil - \lfloor \frac{q-p}{2}\rfloor
 +1$, or 
 $n-\lfloor\frac{n-p}{2}\rfloor+\lceil \frac{q-p}{2}\rceil +1
 \le u\le n$. 

Analogously, a symmetric $n\times n$ matrix $A^{\prime}$ is called a {\sl permuted-shifted
monotone Monge matrix} ({\sl PS monotone Monge matrix }), if it can be obtained form a
permuted-shifted monotone anti-Monge $A$ by multiplying it by $-1$ and by then
adding a sum matrix to it (which can also be the
zero matrix, i.e.\ the matrix containig only entries equal to zero). 
See Figure~\ref{fig:8Rays} for a graphical illustration of PS monotone
anti-Monge  and PS monotone Monge matrices.
\end{definition}
\smallskip

\noindent Summarizing we have proved the following result:
\begin{theorem}
The $QAP(A,B)$ with a PS monotone anti-Monge matrix $A$ and an up-benevolent Toeplitz matrix $B$ is
solved to optimality by the identity permutation.
\end{theorem}

Notice that this result is a strict  generalisation of the result of Burkard \&
al.~\cite{BCRW1998}, because if a PS-anti-Monge matrix A is permuted by
the inverse $(\pi^*)^{-1}$ of the Supnick permutation $\pi^*$, it generally  does yield an
anti-Monge matrix. 

Finally observe that, clearly,  the $n\times n$ PS monotone anti-Monge matrices  form also
cone whose extremal rays  are the matrices $C^{(p,q,u)}$  with $1\le p\le q\le n$,  $1< u\le \lceil \frac{n-p}{2}\rceil - \lfloor \frac{q-p}{2}\rfloor
 +1$, or 
 $n-\lfloor\frac{n-p}{2}\rfloor+\lceil \frac{q-p}{2}\rceil +1
 \le u\le n$.
 Thus these extremal rays build a three
parametric family of matrices  in contrast to the extremal rays of the
(symmetric) monotone anti-Monge matrices which build a two paramteric family.

Since the equality  $Z(A,B,\pi)=Z(-A,-B,\pi)$ trivially holds for any
permutation $\pi$,  $QAP(A,B)$ and $QAP(-A,-B)$ have the same set of optimal
solutions. Notice,  moreover, that if $B$ is up-benevolent Toeplitz matrix than
$-B$ is a down-benevolent Toeplitz matrix. Summarizing we obtain

\begin{corollary}\label{PS-monge:corro}
The $QAP(A,B)$ with a PS monotone Monge matrix $A$ and a down-benevolent Toeplitz  matrix $B$ is
solved to optimality by the identity permutation.
\end{corollary}

\subsection{Combined QAPs}\label{combined:ssec}

In the previous sections we reviewed known  polynomially solvable cases of the
QAP and  proved some new results. 
In this section we show that some of the special  cases use the same special
structures, and can hence  be combined into new structures and new solvable cases. 

\begin{figure}[htb]
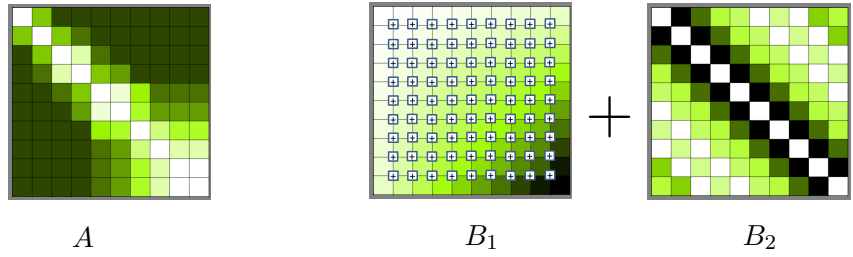

\unitlength=1cm
\begin{center}
\begin{picture}(11.5,2.5)
\put(0,0.5){\begin{picture}(3.8,3.)
\includegraphics[scale=0.5]{BlockCombined}
\end{picture} }
\put(1,0){\makebox(0,0)[cc]{$A$}}
\put(8,1.6){\makebox(0,0)[cc]{\Huge +}}
\put(4.8,0.5){\begin{picture}(4.2,3.)
\includegraphics[scale=0.5]{AntiMonge}
\end{picture} }
\put(6.3,0){\makebox(0,0)[cc]{$B_1$}}
\put(8.5,0.5){\begin{picture}(4.2,3.)
\includegraphics[scale=0.5]{Malevolent}
\end{picture} }
\put(10,0){\makebox(0,0)[cc]{$B_2$}}
\end{picture}
\end{center}
\caption{Illustration to $QAP(A,B_1+B_2)$ where $A$ is a conic combination of cut matrices in CDW normal form, $B_1$ is a monotone anti-Monge matrix, $B_2$ is a down-benevolent Toeplitz matrix.}
\label{fig:Combined1}
\end{figure}

\smallskip
{\bf Cut matrices in CDW normal form.} Consider a $QAP(A,B)$ where  the matrix $A$
is  a conic combination of cut matrices in CDW normal form. Since this matrix
is both a  Kalmanson and a Robinson matrix, we can combine the special  case
described in \cite{Cela2015} and the new special presented in
Section~\ref{down-benev:ssec} matrix $B$ can now be chosen to be  a conic combination of two matrices - a
 symmetric monotone  anti-Monge  matrix and a down-benevolent Toeplitz
matrix 
(see illustration on Figure~\ref{fig:Combined1}).
\bigskip
\bigskip

\begin{figure}[htb]
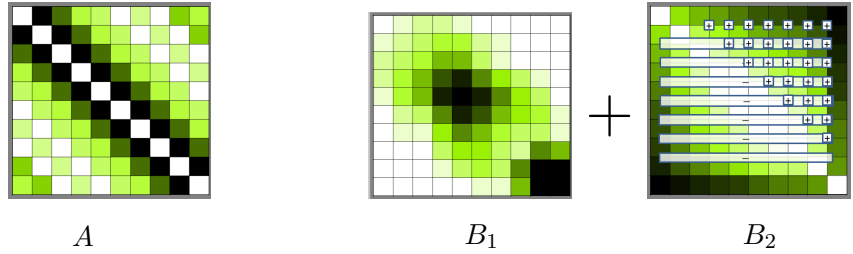

\unitlength=1cm
\begin{center}
\begin{picture}(11.5,3.5)
\put(0,0.5){\begin{picture}(3.8,4)
\includegraphics[scale=0.5]{Malevolent}
\end{picture} }
\put(1,0){\makebox(0,0)[cc]{$A$}}
\put(8,1.6){\makebox(0,0)[cc]{\Huge +}}
\put(4.8,0.5){\begin{picture}(4.2,3.6)
\includegraphics[scale=0.5]{PSantiM}
\end{picture} }
\put(6.3,0){\makebox(0,0)[cc]{$B_1$}}
\put(8.5,0.5){\begin{picture}(4.2,3.6)
\includegraphics[scale=0.5]{KalmSorted}
\end{picture} }
\put(10,0){\makebox(0,0)[cc]{$B_2$}}
\end{picture}
\end{center}
\caption{Illustration to $QAP(A,B_1+B_2)$ where $A$ is a down-benevolent
  Toeplitz matrix, $B_1$ is a PS monotone Monge matrix, $B_2$ is a Kalmanson and Robinsonian dissimilarity matrix.}
\label{fig:Combined2}
\end{figure}

{\bf Down-benevolent  Toeplitz.}
 By combining the special cases described in   Section~\ref{down-benev:ssec} 
 and in Corollary~\ref{PS-monge:corro}  we get a new solvable case of the $QAP(A,B_1+B_2)$ where
$A$ is a down-benevolent  Toeplitz matrix, and the second matrix $B$ is a conic
combination of  matrices wich are both Kalmanson and Robinson matrices  and a
PS monotone Monge matrix  (see illustration on Figure~\ref{fig:Combined2}).

\medskip

{\bf  DW-Toeplitz.}
Let $A$ be an $n\times n$ DW-Toeplitz matrix with generating function $f$. By definition $A$
is also a circulant matrix and
fulfills $f(i)=f(n-i)=f(i-n)$, for $i=,1,2,\ldots,n-1$. Clearly such a matrix
$A$ is a special down-benevolent matrix. 
 Consider now a  special case $QAP(A,B)$ where the identity is an
 optimal solution. 
Since matrix $A$ has a circular structure, the identity is still an optimal
solution of  $QAP(A,B^{(u)})$, where $B^{(u)}$ is obtained from $B$ by
applying to it an arbitrary  cyclic shift  according to some permutation $\sigma_u=\langle
   u,u+1,\ldots,n,1,\ldots,u-1  \rangle$,  for any $1\le u \le n $ ($u=1$
   yields  the
   identity permutation as a trivial cyclic shift). 
 In particular consider a $QAP(A,B)$, where $A$ is a DW-Toeplitz matrix and
 $B=-\bar{R}^{(p,q)}$,
 for some $1\le p\le q\le n$. This QAP is solved to optimality by the identity
 permutation as mentioned in Section~\ref{up-benev:ssec}. 
But then the identity permutation is also  an optimal solution of the
$QAP(A,-C^{(p,q,u)})$, where $-C^{(p,q,u)}$ is obtained from
$B=-\bar{R}^{(p,q)}$ by permuting it according to $\sigma_u$, $1\le u\le n$. 
Thus we can extend the
   class of PS monotone Monge  matrices
 defined  in Section~\ref{up-benev:ssec} and obtain the class of  the {\sl cyclic PS monotone Monge
   matrices},  which is  the  class of matrices
obtained  by first permuting $-R^{(p,q)}$ according to
$\pi^{\ast}$ and then by permuting the resulting matrix according to a cyclic shift $\sigma_u$,   $1\le p\le
q\le n$ and $1\le u\le n$.

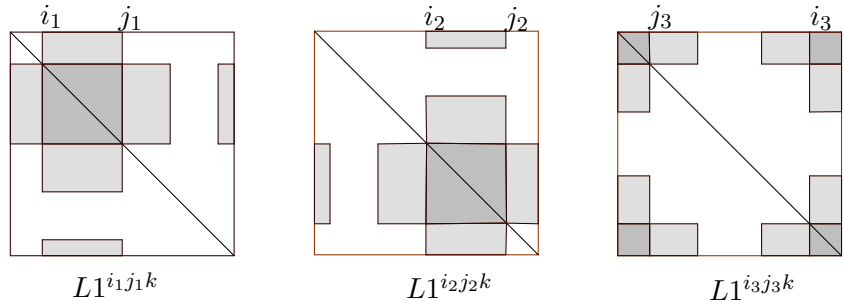
\begin{figure}[htb]
\unitlength=1cm
\begin{center}
\begin{picture}(11,3.65)
\put(-0.3,0)
{
\definecolor{zzttqq}{rgb}{0.6,0.2,0.}
\definecolor{ttqqqq}{rgb}{0.2,0.,0.}
\begin{tikzpicture}[line cap=round,line join=round,>=triangle 45,x=0.85cm,y=0.85cm]
\clip(-2.,8.5) rectangle (12.,14.);
\fill[fill=gray!50] (-1.,11.25) -- (0.25,11.25) -- (0.25,12.5) -- (-1.,12.5) -- cycle;
\fill[fill=gray!25] (-1.,12.5) -- (-1.,13.) -- (0.25,13.) -- (0.25,12.5) -- cycle;
\fill[fill=gray!25] (-1.5,12.5) -- (-1.,12.5) -- (-1.,11.25) -- (-1.5,11.25) -- cycle;
\fill[fill=gray!25] (-1.,11.25) -- (-1.,10.5) -- (0.25,10.5) -- (0.25,11.25) -- cycle;
\fill[fill=gray!25] (0.25,12.5) -- (1.,12.5) -- (1.,11.25) -- (0.25,11.25) -- cycle;
\fill[fill=gray!25] (-1.,9.75) -- (0.25,9.75) -- (0.25,9.5) -- (-1.,9.5) -- cycle;
\fill[fill=gray!25] (1.75,12.5) -- (2.,12.5) -- (2.,11.25) -- (1.75,11.25) -- cycle;
\fill[fill=gray!50] (6.2604790780560196,10.010479078056022) -- (5.,10.) -- (5.01047907805602,11.260479078056022) -- (6.25,11.25) -- cycle;
\fill[fill=gray!25] (5.01047907805602,11.260479078056022) -- (5.,12.) -- (6.25,12.) -- (6.25,11.25) -- cycle;
\fill[fill=gray!25] (4.25,11.25) -- (5.01047907805602,11.260479078056022) -- (5.,10.) -- (4.25,10.) -- cycle;
\fill[fill=gray!25] (6.25,11.25) -- (6.756583579114702,11.25) -- (6.756583579114702,10.) -- (6.2604790780560196,10.010479078056022) -- cycle;
\fill[fill=gray!25] (5.,10.) -- (6.2604790780560196,10.010479078056022) -- (6.25,9.51437457699734) -- (5.,9.51437457699734) -- cycle;
\fill[fill=gray!25] (3.256583579114699,11.25) -- (3.5,11.25) -- (3.5,10.) -- (3.256583579114699,10.) -- cycle;
\fill[fill=gray!25] (5.,13.014374576997344) -- (6.25,13.014374576997344) -- (6.25,12.75) -- (5.,12.75) -- cycle;
\fill[fill=gray!50] (11.,10.) -- (11.5,10.) -- (11.5,9.5) -- (11.,9.5) -- cycle;
\fill[fill=gray!50] (11.,13.) -- (11.,12.5) -- (11.5,12.5) -- (11.5,13.) -- cycle;
\fill[fill=gray!50] (8.,13.) -- (8.493364613070728,13.) -- (8.5,12.5) -- (8.,12.5) -- cycle;
\fill[fill=gray!50] (8.,10.) -- (8.5,10.) -- (8.5,9.5) -- (8.,9.5) -- cycle;
\fill[fill=gray!25] (11.,10.) -- (11.,10.75) -- (11.5,10.75) -- (11.5,10.) -- cycle;
\fill[fill=gray!25] (11.,12.5) -- (11.,11.75) -- (11.5,11.754821485011426) -- (11.5,12.5) -- cycle;
\fill[fill=gray!25] (8.,12.5) -- (8.5,12.5) -- (8.5,11.75) -- (8.,11.75) -- cycle;
\fill[fill=gray!25] (8.5,12.5) -- (9.25,12.5) -- (9.25,13.) -- (8.493364613070728,13.) -- cycle;
\fill[fill=gray!25] (8.,10.) -- (8.,10.75) -- (8.5,10.75) -- (8.5,10.) -- cycle;
\fill[fill=gray!25] (8.5,10.) -- (9.25,10.) -- (9.25,9.5) -- (8.5,9.5) -- cycle;
\fill[fill=gray!25] (11.,9.5) -- (10.25,9.5) -- (10.25,10.) -- (11.,10.) -- cycle;
\fill[fill=gray!25] (10.25,13.) -- (11.,13.) -- (11.,12.5) -- (10.25,12.5) -- cycle;
\draw [color=ttqqqq] (-1.5,9.5)-- (2.,9.5);
\draw [color=ttqqqq] (2.,9.5)-- (2.,13.);
\draw [color=ttqqqq] (2.,13.)-- (-1.5,13.);
\draw [color=ttqqqq] (-1.5,13.)-- (-1.5,9.5);
\draw (-1.5,13.)-- (2.,9.5);
\draw (-1.19171881802166,13.585371507864773) node[anchor=north west] {$i_1$};
\draw (0.03476201821527471,13.566788464891486) node[anchor=north west] {$j_1$};
\draw (-0.689976657742914,9.4041868388752) node[anchor=north west] {$L1^{i_1j_1k}$};
\draw [color=zzttqq] (3.256583579114699,9.51437457699734)-- (6.756583579114702,9.51437457699734);
\draw [color=zzttqq] (6.756583579114702,9.51437457699734)-- (6.756583579114702,13.014374576997344);
\draw [color=zzttqq] (6.756583579114702,13.014374576997344)-- (3.2565835791146984,13.014374576997344);
\draw [color=zzttqq] (3.2565835791146984,13.014374576997344)-- (3.256583579114699,9.51437457699734);
\draw (3.2565835791146984,13.014374576997344)-- (6.756583579114702,9.51437457699734);
\draw (4.810604062350005,13.566788464891486) node[anchor=north west] {$i_2$};
\draw (6.0742509845335135,13.548205421918198) node[anchor=north west] {$j_2$};
\draw (10.850093028668244,13.566788464891486) node[anchor=north west] {$i_3$};
\draw (8.322799184301227,13.585371507864773) node[anchor=north west] {$j_3$};
\draw (4.531858417750702,9.385603795901913) node[anchor=north west] {$L1^{i_2j_2k}$};
\draw [color=ttqqqq] (-1.,11.25)-- (0.25,11.25);
\draw [color=ttqqqq] (0.25,11.25)-- (0.25,12.5);
\draw [color=ttqqqq] (0.25,12.5)-- (-1.,12.5);
\draw [color=ttqqqq] (-1.,12.5)-- (-1.,11.25);
\draw [color=ttqqqq] (-1.,12.5)-- (-1.,13.);
\draw [color=ttqqqq] (-1.,13.)-- (0.25,13.);
\draw [color=ttqqqq] (0.25,13.)-- (0.25,12.5);
\draw [color=ttqqqq] (0.25,12.5)-- (-1.,12.5);
\draw [color=ttqqqq] (-1.5,12.5)-- (-1.,12.5);
\draw [color=ttqqqq] (-1.,12.5)-- (-1.,11.25);
\draw [color=ttqqqq] (-1.,11.25)-- (-1.5,11.25);
\draw [color=ttqqqq] (-1.5,11.25)-- (-1.5,12.5);
\draw [color=ttqqqq] (-1.,11.25)-- (-1.,10.5);
\draw [color=ttqqqq] (-1.,10.5)-- (0.25,10.5);
\draw [color=ttqqqq] (0.25,10.5)-- (0.25,11.25);
\draw [color=ttqqqq] (0.25,11.25)-- (-1.,11.25);
\draw [color=ttqqqq] (0.25,12.5)-- (1.,12.5);
\draw [color=ttqqqq] (1.,12.5)-- (1.,11.25);
\draw [color=ttqqqq] (1.,11.25)-- (0.25,11.25);
\draw [color=ttqqqq] (0.25,11.25)-- (0.25,12.5);
\draw [color=ttqqqq] (-1.,9.75)-- (0.25,9.75);
\draw [color=ttqqqq] (0.25,9.75)-- (0.25,9.5);
\draw [color=ttqqqq] (0.25,9.5)-- (-1.,9.5);
\draw [color=ttqqqq] (-1.,9.5)-- (-1.,9.75);
\draw [color=ttqqqq] (1.75,12.5)-- (2.,12.5);
\draw [color=ttqqqq] (2.,12.5)-- (2.,11.25);
\draw [color=ttqqqq] (2.,11.25)-- (1.75,11.25);
\draw [color=ttqqqq] (1.75,11.25)-- (1.75,12.5);
\draw [color=ttqqqq] (6.2604790780560196,10.010479078056022)-- (5.,10.);
\draw [color=ttqqqq] (5.,10.)-- (5.01047907805602,11.260479078056022);
\draw [color=ttqqqq] (5.01047907805602,11.260479078056022)-- (6.25,11.25);
\draw [color=ttqqqq] (6.25,11.25)-- (6.2604790780560196,10.010479078056022);
\draw [color=ttqqqq] (5.01047907805602,11.260479078056022)-- (5.,12.);
\draw [color=ttqqqq] (5.,12.)-- (6.25,12.);
\draw [color=ttqqqq] (6.25,12.)-- (6.25,11.25);
\draw [color=ttqqqq] (6.25,11.25)-- (5.01047907805602,11.260479078056022);
\draw [color=ttqqqq] (4.25,11.25)-- (5.01047907805602,11.260479078056022);
\draw [color=ttqqqq] (5.01047907805602,11.260479078056022)-- (5.,10.);
\draw [color=ttqqqq] (5.,10.)-- (4.25,10.);
\draw [color=ttqqqq] (4.25,10.)-- (4.25,11.25);
\draw [color=ttqqqq] (6.25,11.25)-- (6.756583579114702,11.25);
\draw [color=ttqqqq] (6.756583579114702,11.25)-- (6.756583579114702,10.);
\draw [color=ttqqqq] (6.756583579114702,10.)-- (6.2604790780560196,10.010479078056022);
\draw [color=ttqqqq] (6.2604790780560196,10.010479078056022)-- (6.25,11.25);
\draw [color=ttqqqq] (5.,10.)-- (6.2604790780560196,10.010479078056022);
\draw [color=ttqqqq] (6.2604790780560196,10.010479078056022)-- (6.25,9.51437457699734);
\draw [color=ttqqqq] (6.25,9.51437457699734)-- (5.,9.51437457699734);
\draw [color=ttqqqq] (5.,9.51437457699734)-- (5.,10.);
\draw [color=ttqqqq] (3.256583579114699,11.25)-- (3.5,11.25);
\draw [color=ttqqqq] (3.5,11.25)-- (3.5,10.);
\draw [color=ttqqqq] (3.5,10.)-- (3.256583579114699,10.);
\draw [color=ttqqqq] (3.256583579114699,10.)-- (3.256583579114699,11.25);
\draw [color=ttqqqq] (5.,13.014374576997344)-- (6.25,13.014374576997344);
\draw [color=ttqqqq] (6.25,13.014374576997344)-- (6.25,12.75);
\draw [color=ttqqqq] (6.25,12.75)-- (5.,12.75);
\draw [color=ttqqqq] (5.,12.75)-- (5.,13.014374576997344);
\draw [color=zzttqq] (8.,13.)-- (8.,9.5);
\draw [color=zzttqq] (8.,9.5)-- (11.5,9.5);
\draw [color=zzttqq] (11.5,9.5)-- (11.5,13.);
\draw [color=zzttqq] (11.5,13.)-- (8.,13.);
\draw (8.,13.)-- (11.5,9.5);
\draw [color=ttqqqq] (11.,10.)-- (11.5,10.);
\draw [color=ttqqqq] (11.5,10.)-- (11.5,9.5);
\draw [color=ttqqqq] (11.5,9.5)-- (11.,9.5);
\draw [color=ttqqqq] (11.,9.5)-- (11.,10.);
\draw [color=ttqqqq] (11.,13.)-- (11.,12.5);
\draw [color=ttqqqq] (11.,12.5)-- (11.5,12.5);
\draw [color=ttqqqq] (11.5,12.5)-- (11.5,13.);
\draw [color=ttqqqq] (11.5,13.)-- (11.,13.);
\draw [color=ttqqqq] (8.,13.)-- (8.493364613070728,13.);
\draw [color=ttqqqq] (8.493364613070728,13.)-- (8.5,12.5);
\draw [color=ttqqqq] (8.5,12.5)-- (8.,12.5);
\draw [color=ttqqqq] (8.,12.5)-- (8.,13.);
\draw [color=ttqqqq] (8.,10.)-- (8.5,10.);
\draw [color=ttqqqq] (8.5,10.)-- (8.5,9.5);
\draw [color=ttqqqq] (8.5,9.5)-- (8.,9.5);
\draw [color=ttqqqq] (8.,9.5)-- (8.,10.);
\draw [color=ttqqqq] (11.,10.)-- (11.,10.75);
\draw [color=ttqqqq] (11.,10.75)-- (11.5,10.75);
\draw [color=ttqqqq] (11.5,10.75)-- (11.5,10.);
\draw [color=ttqqqq] (11.5,10.)-- (11.,10.);
\draw [color=ttqqqq] (11.,12.5)-- (11.,11.75);
\draw [color=ttqqqq] (11.,11.75)-- (11.5,11.754821485011426);
\draw [color=ttqqqq] (11.5,11.754821485011426)-- (11.5,12.5);
\draw [color=ttqqqq] (11.5,12.5)-- (11.,12.5);
\draw [color=ttqqqq] (8.,12.5)-- (8.5,12.5);
\draw [color=ttqqqq] (8.5,12.5)-- (8.5,11.75);
\draw [color=ttqqqq] (8.5,11.75)-- (8.,11.75);
\draw [color=ttqqqq] (8.,11.75)-- (8.,12.5);
\draw [color=ttqqqq] (8.5,12.5)-- (9.25,12.5);
\draw [color=ttqqqq] (9.25,12.5)-- (9.25,13.);
\draw [color=ttqqqq] (9.25,13.)-- (8.493364613070728,13.);
\draw [color=ttqqqq] (8.493364613070728,13.)-- (8.5,12.5);
\draw [color=ttqqqq] (8.,10.)-- (8.,10.75);
\draw [color=ttqqqq] (8.,10.75)-- (8.5,10.75);
\draw [color=ttqqqq] (8.5,10.75)-- (8.5,10.);
\draw [color=ttqqqq] (8.5,10.)-- (8.,10.);
\draw [color=ttqqqq] (8.5,10.)-- (9.25,10.);
\draw [color=ttqqqq] (9.25,10.)-- (9.25,9.5);
\draw [color=ttqqqq] (9.25,9.5)-- (8.5,9.5);
\draw [color=ttqqqq] (8.5,9.5)-- (8.5,10.);
\draw [color=ttqqqq] (11.,9.5)-- (10.25,9.5);
\draw [color=ttqqqq] (10.25,9.5)-- (10.25,10.);
\draw [color=ttqqqq] (10.25,10.)-- (11.,10.);
\draw [color=ttqqqq] (11.,10.)-- (11.,9.5);
\draw [color=ttqqqq] (10.25,13.)-- (11.,13.);
\draw [color=ttqqqq] (11.,13.)-- (11.,12.5);
\draw [color=ttqqqq] (11.,12.5)-- (10.25,12.5);
\draw [color=ttqqqq] (10.25,12.5)-- (10.25,13.);
\draw (9.270534375938858,9.367020752928626) node[anchor=north west] {$L1^{i_3j_3k}$};
\end{tikzpicture}
}
\end{picture}
\end{center}
\caption{Illustration to rays which generate Circular Shift and Permuted Monge matrix.}
\label{fig:LRaysCircle}
\end{figure}

 Figure \ref{fig:LRaysCircle} illustrates such cyclic PS monotone Monge matrices.
\smallskip

\begin{figure}[htb]
\unitlength=1cm
\begin{center}
\begin{picture}(11.5,3.5)
\put(0,0.5){\begin{picture}(3.8,3.5)
\includegraphics[scale=0.5]{ToeplitzDW}
\end{picture} }
\put(1,0){\makebox(0,0)[cc]{$A$}}
\put(8,1.6){\makebox(0,0)[cc]{\Huge +}}
\put(4.8,0.5){\begin{picture}(4.2,3.5)
\includegraphics[scale=0.5]{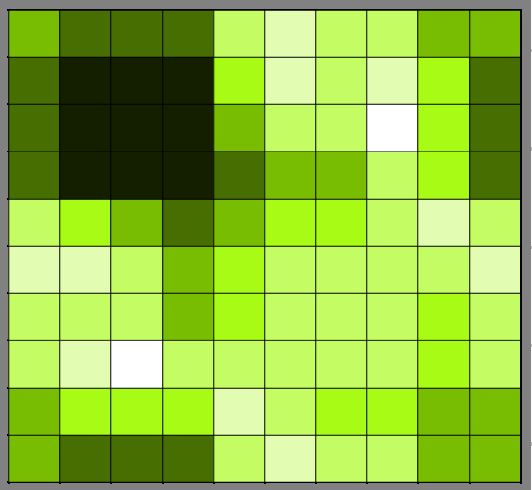}
\end{picture} }
\put(6.3,0){\makebox(0,0)[cc]{$B_1$}}
\put(8.5,0.5){\begin{picture}(4.2,3.5)
\includegraphics[scale=0.5]{KalmansonUnsorted}
\end{picture} }
\put(10,0){\makebox(0,0)[cc]{$B_2$}}
\end{picture}
\end{center}
\caption{Illustration to $QAP(A,B_1+B_2)$ where $A$ is a DW-Toeplitz matrix,
  $B_1$ is a cyclic PS monotone Monge matrix, $B_2$ is a Kalmanson matrix.}
\label{fig:Combined3}
\end{figure}

We can define now a new combined special case of the QAP  solved to
optimality by the identity permutation, namely  $QAP(A,B_1+B_2)$,  where $A$ is
a DW-Toeplitz matrix, and $B$ is a conic combination
of a Kalmanson matrix  and a cyclic PS monotone Monge matrix (see the illustration in Figure~\ref{fig:Combined3}).

\bigskip
\smallskip

\section{Conic representation of specially structured matrices}
\label{sec:cut} 
\subsection{Cut weights and specially structured matrices}\label{cutweights:ssec}
In this section, we investigate the structure of matrices which are both Kalmanson  and Robinson matrices. 
We show that any matrix in this class can be represented as a sum of a constant
matrix and a conic combination of 
cut matrices.


We use Lemma~\ref{lemma: kalm.1} as an alternative definition  of Kalmason matrices.
Consider special  cut matrices  $A^{kl}$, $1\le k<l\le n$,  containing  one block of size $(k-l+1)$ with $a_{ij}=0$ for $k\le i,j\le l$, and all other $n-k+l-1$ blocks of size $1$.

Is can be easily observed that  the matrices $A^{(k,l)}$ fulfill the
inequalties \ref{kalmx.c1} and \ref{kalmx.c1} and are therefore Kalmanson
matrices. 
Notice moreover  that for any $n\times n$ cut matrix $A^{kl}$, $1<k<l<n$, there
is only one strict inequality 
in \eqref{kalmx.c1}, namely 
\[ a_{k-1,l}+a_{k,l+1} > a_{kl}+a_{k-1,l+1}\, , \]
whereas all inequalities \eqref{kalmx.c2}
 are fulfilled with equality.
Analogously, there is only one strict inequality in \eqref{kalmx.c2} for the
matrices  $A^{1,k-1}$ and $A^{kn}$, $2< k<n$, namely
\[ a_{k-1,1}+a_{kn} < a_{k1}+a_{k-1,n}\, , \]
whereas all  inequalities \eqref{kalmx.c1}  are fulfilled with equality.

The following lemma shows that any  Kalmanson matrix  can be represented as a
linear combination of a weak sum matrix with cut matrices $A^{k,l}$ which can be computed explicitly
in terms of  simple formulas involving the entries of the considered Kalmanson matrix. 
Similar structural properties of Kalmanson matrices in terms of cuts and cut-weights have also been 
studied in \cite{BaDr} and \cite{ChFaTr}.
In both papers though the authors suggest algorithms for calculating the cut-weights while we 
provide simple analytical expressions for them.
\smallskip

\begin{lemma}
\label{lemma:kalm.2}
A symmetric $n\times n$ matrix $C$ is a Kalmanson matrix if and only if it can be represented as a 
linear combination of a weak sum matrix $S$ and cut matrices $A^{k,l}$ as follows 
\begin{eqnarray}\label{kalm2}
C=S+\sum_{i=1}^{n-3}\sum_{j=i+2}^{n-1}
 \delta_{i+1,j} A^{i+1,j}
+ \sum_{i=2}^{n-2}(\alpha_i A^{1,i}+\beta_i A^{i+1,n})\, .
\end{eqnarray}
The coefficients of the linear combination, the so-called {\sl cut weights},
are given as 
$\delta_{i+1,j}=(-c_{i,j+1}-c_{i+1,j} + c_{ij}+c_{i+1,j+1})$,
$\alpha_i=c_{i+1,1}-c_{i,1},\ \beta_i=c_{in}-c_{i+1,n}$ and fulfill $\delta_{i+1,j}\ge 0$, $\alpha_i+\beta_i\ge 0$.
\end{lemma}

\proof 
It can easily be checked that any weak sum matrix, a cut matrix $A^{kl}$, and a linear combination 
$\alpha_i A^{1,i}+\beta_i A^{i+1,n}$ with $\alpha_i+\beta_i\ge 0$ are Kalmanson matrices, and therefore 
any matrix given  as  in \eqref{kalm2} is a Kalmanson matrix.

Assume now that $C$ is a Kalmanson matrix. Let $i$ and $j$, $1\le i<i+2\le
j<n-1$ be two indices, such the corresponding inequality  in \eqref{kalmx.c1}
is strict, i.e.\
$c_{i,j+1}+c_{i+1,j} <c_{ij}+c_{i+1,j+1}$ holds. The involved matrix entries
are printed in boldface in the illustration below, note that all these entries
lie above the main diagonal.

\[
C=\left( \begin{array}{ccccccc}
&&&\ldots &&&\\
 \ldots\  & c_{i,p}&\ldots &c_{i,j-1}&  \mathbf c_{i,j}&\mathbf c_{i,j+1}&\ldots\\
 \ldots\ & c_{i+1,p}&\ldots &c_{i+1,j-1}&\mathbf c_{i+1,j}&\mathbf c_{i+1,j+1}&\ldots\\ 
 &&&\ldots &&&\\
  \ldots\ & c_{q,p}&\ldots &c_{q,j-1}&c_{q,j}&c_{q,j+1}&\ldots\\ 
&&&\ldots &&&\\
\end{array} \right)
\]
Set $\delta_{i+1,j}:=-c_{i,j+1}-c_{i+1,j} + c_{ij}+c_{i+1,j+1}>0$ and consider
the  matrix $C'=C-\delta_{i+1,j}A^{i+1,j}$, represented schematically below (to
simplify the illustration 
we use the notation $\Delta:=\delta_{i+1,j}$): 
\smallskip

\[
C'=\left( \begin{array}{ccccccccc}
c_{11}-\Delta& \ldots & c_{1,i}-\Delta & c_{1,i+1}-\Delta & \ldots &
c_{1j}-\Delta& c_{1,j+1}-\Delta&\ldots & c_{1,n}-\Delta\\
\ddots & \ddots & \ddots & \ddots  & \ddots &
\ddots  & \ddots &\ddots & \ddots  \\
c_{i,1}-\Delta& \ldots & c_{i,i}-\Delta & c_{i,i+1}-\Delta & \ldots &
\mathbf c_{i,j}-\Delta& \mathbf c_{i,j+1}-\Delta&\ldots & c_{i,n}-\Delta\\
c_{i+1,1}-\Delta& \ldots & c_{i+1,i}-\Delta & c_{i+1,i+1} & \ldots &
\mathbf c_{i+1,j}& \mathbf c_{i+1,j+1}-\Delta&\ldots & c_{i+1,n}-\Delta\\
\ddots & \ddots & \ddots & \ddots  & \ddots &
\ddots  & \ddots &\ddots & \ddots  \\
c_{j,1}-\Delta& \ldots & c_{j,i}-\Delta & c_{j,i+1} & \ldots &
c_{j,j}& c_{j,j+1}-\Delta&\ldots & c_{j,n}-\Delta\\
c_{j+1,1}-\Delta& \ldots & c_{j+1,i}-\Delta & c_{j+1,i+1}- \Delta & \ldots &
c_{j+1,j}-\Delta& c_{j+1,j+1}-\Delta&\ldots & c_{j+1,n}-\Delta\\
\ddots & \ddots & \ddots & \ddots  & \ddots &
\ddots  & \ddots &\ddots & \ddots  \\
c_{n,1}-\Delta& \ldots & c_{n,i}-\Delta & c_{n,i+1}- \Delta & \ldots &
c_{n,j}-\Delta& c_{n,j+1}-\Delta&\ldots & c_{n,n}-\Delta\\
\end{array} \right)
\] 
\smallskip


Notice  that in matrix $C'=(c'_{ij})$ we have
$c'_{i,j+1}+c'_{i+1,j} =c'_{ij}+c'_{i+1,j+1}$. 
Moreover the status  of the other inequalities in \eqref{kalmx.c1} does not
change, meaning  that all inequalities are still fulfilled by matrix $C'$ and
only the inequalities which were strictly fulfilled by $C$ are strictly
fulfilled by $C'$. Finally it is also easy to see that $C'$ fulfills
inequalities \eqref{kalmx.c1}.  Hence $C'$ is a Kalmanson matrix,  and we check
again whether there is a pair of indices for which the corresponding inequality
in   \eqref{kalmx.c1} is strict. If yes,  we  perfom an analogous transformation
as the one described above by defining  the corresponding $\delta$-coefficient
and substracting from $C'$ the corresponding cut matrix multiplied by that coefficient. We repeat this process, update $C'$ in every step,  and eventually obtain a Kalmanson matrix $C'$
which filfills all inequalities \eqref{kalmx.c1} by equality. 

Assume now that there exists some inequality in  \eqref{kalmx.c2}  strictly
fulfilled by the entries of $C'$. 
Let  $i$, $2\le i \le n-2$, be an index  such that $c_{i1}+c_{i+1,n}<c_{i+1,1}+c_{in}$:
\[
C=\left( \begin{array}{ccccccc}
&&&\ldots &&&\\
\mathbf c_{i1}&c_{i2}&\ldots & 0\ \  &\ldots &  c_{i,n-1}&\mathbf c_{in}\\
\mathbf c_{i+1,1}&c_{i+1,2}&\ldots &\ \  0 &\ldots &  c_{i+1,n-1}&\mathbf c_{i+1,n}\\
&&&\ldots &&&\\
\end{array} \right)
\]

Set $\alpha_i=c_{i+1,1}-c_{i,1}$, $\beta_i=c_{in}-c_{i+1,n}$. 
 Clearly  $\alpha_i+\beta_i>0$ holds, due to
 $c_{i1}+c_{i+1,n}<c_{i+1,1}+c_{in}$.   
Consider the matrix  $C'=C-\alpha A^{1i}-\beta\ A^{i+1,n}$ where
$\alpha:=\alpha_i$, $\beta:=\beta_i$:
\smallskip

\[ 
C'=\left( \begin{array}{cccccc}
c_{1,1}-\beta & \ldots & c_{1,i}-\beta & c_{1,i+1}-\alpha-\beta & \ldots &
  c_{1,n}-\alpha-\beta \\
\ddots & \ddots &\ddots &\ddots &\ddots &\ddots \\
\mathbf c_{i,1}-\beta &\ldots &c_{i,i}-\beta & c_{i,i+1}-\alpha-\beta &\ldots &
 \mathbf c_{i,n}-\alpha-\beta \\
\mathbf c_{i+1,1}-\alpha -\beta &\ldots &c_{i+1,i}-\alpha - \beta & c_{i+1,i+1}-\alpha &\ldots &
  \mathbf c_{i+1,n}-\alpha \\
\ddots & \ddots &\ddots &\ddots &\ddots &\ddots \\
c_{n,1}-\alpha -\beta &\ldots &c_{n,i}-\alpha - \beta & c_{n,i+1}-\alpha &\ldots &
  c_{n,n}-\alpha \\
\end{array} \right)
\]
\smallskip

It can be easily checked that $c'_{i1}+c'_{i+1,n}=c'_{i+1,1}+c'_{in}$,  that
all inequalities \eqref{kalmx.c1} remain fulfilled with  equality,  and that
the status   of the other inequalities in \eqref{kalmx.c2} does not
change,  meaning  that all these  inequalities are still fulfilled by matrix $C'$ and
only those inequalities among them which were strictly fulfilled by $C$, are strictly
fulfilled by $C'$. As long as there are inequalities \eqref{kalmx.c2} strictly
fulfilled by $C'$ we apply a transformation as above on $C'$ and update $C'$. 
So eventually we get a transformed matrix where all inequalities 
\eqref{kalmx.c1}, \eqref{kalmx.c2} are fulfilled with equality. 
Such a matrix is a weak sum matrix, as  shown in Lemma~\ref{weaksum:lemma} below, and
this completes the proof.
\qed
\smallskip
\begin{lemma}\label{weaksum:lemma}
Let $C$ be an $n\times n $ Kalmanson matrix for which all inequalities in
\eqref{kalmx.c1} and \eqref{kalmx.c2} are fulfilled with equality. Then $C$ is
a weak sum matrix. 
\end{lemma}
\proof
Consider an index $i$, $1\le i \le n-1$. Since \eqref{kalmx.c1} and
\eqref{kalmx.c2} are fulfilled with equality the differences
$c_{i+1,j}-c_{i,j}=c_{j,i+1}-c_{j,i}$,  have  a common value for all $j\in
\{1,2,\ldots,n\}\setminus \{i,i+1\}$. Denote this common value by $b_i$,
$1\le i\le n-1$. 

\noindent Consider now an entry $c_{ij}$ of $C$ with $1\le i<j$. If $j\ge i+2$, the following
equalities hold
\[ c_{i,j}=c_{i,j-1}+b_{j-1}=\ldots= c_{i,i+1}+b_{i+1}+\ldots +
  b_{j-1}\, . \]
The later equality can be also seen as fulfilled for $j=i+1$ where the sum of the
$b$-s would disappear.
Further if $i>1$ we have 
\[c_{ij}=c_{i-1,i+1}+b_{i-1}+(b_{i+1}+\ldots+b_{j-1})=\ldots=
  c_{1,i+1}+(b_1+\ldots,+b_{i-1})+
  (b_{i+1}+\ldots+b_{j-1})\, .\] 
The last equality can be also seen as fulfilled for $i=1$ where the left-most sum of the
$b$-s would disappear.Further, if $i+1>2$, we get  
\[c_{ij}=c_{1,i}+b_i+(b_1+\ldots,+b_{i-1})+
  (b_{i+1}+\ldots+b_{j-1})=\ldots=\]
\[
c_{12}+(b_2+\ldots,b_i) +   (b_1+\ldots,+b_{i-1})+
  (b_{i+1}+\ldots+b_{j-1})\, .\]
The last equality can be also seen as fulfilled for $i=1$ where the two left-most sum of the
$b$-s would  disappear.
Further, for  $1\le i< j\le n$ we get 
\[c_{ij} =c_{12}-b_1+ \sum_{k=1}^{i-1}b_k+
  \sum_{k=1}^{j-1}b_k=\gamma_i+\gamma_j\]
with $\gamma_i:=\sum_{k=1}^{i-1}b_k+(c_{1,2}-b_1)/2$,  for all $i\ge 2$. 
Define now $\gamma_1$ such that $c_{1,2}=\gamma_1+\gamma_2$, that is
$\gamma_1=c_{12}-\gamma_2=
(c_{12}-b_1)/2$. 
So the entries $c_{ij}$ of the symmetric matrix $C$ can be represented as
$c_{ij}=\gamma_i+\gamma_j$, 
for all $1\le i < j\le n$, and hence $C$ is a weak sum matrix. 
\qed

\smallskip

Next we give a characterisation of Kalmanson matrices which are also Robinson
matrices.

\begin{lemma}
\label{lemma:kalm.3}
A symmetric $n\times n$ Kalmanson matrix $C$ is a Robinson matrix if and only if it can be represented 
as a conic combination of a weak constant matrix $Z$ and cut matrices $A^{kl}$ as follows
\begin{eqnarray}\label{kalm3}
C=Z+\sum_{i=1}^{n-3}\sum_{j=i+2}^{n-1}
 \delta_{i+1,j} A^{i+1,j}
+ \sum_{i=2}^{n-1}\alpha_i A^{1,i}+\sum_{i=1}^{n-2}\beta_i A^{i+1,n}\, , 
\end{eqnarray}
where
$
\delta_{i+1,j}:=(-c_{i,j+1}-c_{i+1,j} + c_{ij}+c_{i+1,j+1})$, for $1\le i\le
n-3$, $i+2\le j\le n-1$,
$\alpha_i:=c_{i+1,1}-c_{i,1}$, for $2\le i\le n-1$,  
$\beta_i:=c_{in}-c_{i+1,n}$, for $1\le i\le n-2$, and $\delta_{i+1,j}\ge 0$, $\alpha_i\ge 0$, $\beta_i\ge 0$.
\end{lemma}
\proof 
The proof of the ``if''-part of the lemma is straightforward; just observe that
all matrices in the conic combination are Kalmanson and Robinson matrices and
that a conic combination preserves the Kalmanson and Robinson properties
because all of them are defined in terms of inequalities involving the entries
of the matrix. 

We prove now the  ``only if''-part.
since $C$ is a Kalmanson matrix it has a representation as stated by 
Lemma~\ref{lemma:kalm.2} in \eqref{kalm2}. 
Observe that \eqref{kalm2} and \eqref{kalm3} differ  on the first summand,
which is a weak sum matrix in  \eqref{kalm2} and constant matrix in
\eqref{kalm3}, and on the range of summation for the third and the fourth
summand (combined in one single summand in  \eqref{kalm2}). 
 We go thorugh  the procedure applyied in
Lemma~\ref{lemma:kalm.2} and show  the matrix $C'$ resulting after each
transformation step is again a Robinson matrix.  The non-negativity of the coeffiecients   $\alpha_i$ and
$\beta_i$ in \eqref{kalm2} would then follow   directly form the definition of a Robinson
matrix. 

Consider first a  transformation of the type $C'=C-\Delta A^{i+1,j}$, where 
$\Delta=\delta_{i+1,j}$. 
We claim that $c_{ip}-\Delta \ge c_{i+1,p}$ for all $p=i+2,\ldots,j$.
Let $\Lambda_p=c_{ip}-c_{i+1,p}$, $p=i+2,\ldots,j+1$. \\

Since $C$ is a Robinson matrix, we have $\Lambda_p\ge 0$. Since $C$ is a Kalmanson matrix, we have 
$c_{ip}+c_{i+1,j+1}-c_{i+1,p}-c_{i,j+1}=\Lambda_p -\Lambda_j\ge 0$ and  $\Lambda_p\ge\Lambda_j$. 
Clearly $\Delta= \Lambda_j-\Lambda_{j+1}$ and therefore
$c_{ip}-\Delta - c_{i+1,p}=\Lambda_p-\Lambda_j+\Lambda_{j+1}\ge 0$, which
proves the  claim.
 The claim that 
$c_{q,j+1}-\Delta\ge c_{qj}$ for all $q=i+1,\ldots,j-1$ can be  proved in a similar
way. 
So the new matrix $C'$ is a Robinson dissimilarity. 

Consider now a  transformation of the type $C'=C-\alpha_iA^{1,i} -
\beta_iA^{i+1,n}$, for $2\le i\le n-2$.  
The Kalmanson inequalities
\eqref{kalmx.c1} ensure that $C'$ is a Robinson matrix.
So, what is left to prove is that $Z=S-\alpha_{n-1}A^{1,n-1}-\beta_1 A^{2,n}$ is a weak
constant matrix, where $S$ is the weak 
sum matrix in the  presentation \eqref{kalm2}.

Since  every transformation step results in a Robinson matrix, as shown above, 
the weak sum matrix $S$ resulting after the last transformation in the proof of
Lemma~\ref{lemma:kalm.2} is  a Robinson matrix, too. 
It is easily  seen that a symmetric  weak  sum matrix $S=(s_{ij})$ with
$s_{ij}=\gamma_i+\gamma_j$, for $1\le i <  j \le n$, is a Robinson
dissimilarity,  if and only if 
$\gamma_1\ge \gamma_2= \ldots = \gamma_{n-1}\le \gamma_n$. Indeed
$s_{1j}=\gamma_1+\gamma_j\le s_{1,j+1}=\gamma_1+\gamma_{j+1}$ implies
$\gamma_j\le \gamma_{j+1}$ , for $j\in \{2,3,\ldots,n-1\}$, and
$s_{i-1,n}=\gamma_{i-1}+\gamma_n\ge s_{i,n} =\gamma_i+\gamma_n$, implies
$\gamma_{i-1}\ge \gamma_i$, for $i\in \{2,3,\ldots,n-1\}$. 
 
After the last transformation the equalities $s_{n,1}-s_{n-1,1}=
\gamma_n-\gamma_{n-1}=c_{n,1}-c_{n-1,1}=\alpha_{n-1}$ and
$s_{1,n}-s_{2,n}=\gamma_1-\gamma_2=c_{1,n}-c_{2,n}=\beta_1$ clearly hold. 
Observe finally that 
$$S-(\gamma_n-\gamma_{n-1})\times A^{1,n-1}-(\gamma_1-\gamma_2)\times A^{2,n}=S-\alpha_{n-1}A^{1,n-1}-\beta_1 A^{2,n}$$
is a weak constant matrix (with all non-diagonal elements equal to
$2\gamma_2$),
which completes the proof.
\qed
\smallskip

By applying the above lemma to compute the coefficients of the conic combination
for a cut matrix (which is a Kalmanson and a Robinson matrix) we obtain

\begin{corollary}\label{cor:1}
 Let $C$ be a cut matrix with $m$ blocks such that  $k$ of them ($k\le m$)
 contain more than one element.
Let the corresponding $k$  row and column blocks $I_1,I_2,\ldots,I_k$, $|I_j|>1,\ \forall j$,  of $C$ 
 be given as 
$I_1=\{i_1=1,\ldots,j_1\}$, 
$I_2=\{i_2,\ldots,j_2\}$, \ldots, 
$I_k=\{i_k,\ldots,j_k\}$, where  $i_l\ge j_{l-1}+1$ and $i_l<j_l$, for $1\le
l\le k$. 
Then $C$ can be represented as 
$C=Z+\sum_{l=1}^{l=k} A^{i_l, j_l}$, where $Z=(z_{ij})$ with $z_{ij}=-(k-1)$ for $i \neq j$.
\end{corollary}

\subsection{Recognizing conic combinations of cut matrices in CDW 
normal form}\label{recconcombCDWnormal:ssec}
As mentioned in Section~\ref{combined:ssec} a combined polynomially solvable special case of
the QAP arises if one of the coefficient matrices is a conic combination of cut
matrices in CDW normal form and the other one is a conic combination of a
symmetric anti-Monge matrix and a down-neneolent Toeplitz matrix. Thus,   given a matrix $C$ which is both Kalmanson matrix and Robinson matrix, it is a question of  interest whether the 
matrix can be represented as a conic combination of cut matrices in CDW normal
form (notice that every cut matrix in CDW normal form is both a Robinson and a
Kalmanson matrix but not vice-versa). 
Note that a weakly  constant matrix with zeroes on the diagonal and constant
$K\in \rz$ elsewhere can be obtained
by multiplying with $K$ a special cut matrix in CDW normal form with all blocks of length one.
In order to formulate a simple rule for recognizing this special subclass of
Kalmanson (and Robinson) matrices,  we will associate to every (Kalmanson and
Robinson) matrix   $C$ an $n\times n$ symmetric {\sl cut-weight matrix} $D(C)=(d_{ij})$ with 
$d_{ij}:=\delta_{ij}=c_{i-1,j}+c_{i,j+1}-c_{ij}-c_{i-1,j+1}$ for $2\le i<j\le n-1$, and 
$d_{1i}:=\alpha_i=c_{i+1,1}-c_{i1}$, $i=2,\ldots,n-1$,
$d_{in}:=\beta_{i-1}=c_{i-1,n}-c_{in}$, $i=2,\ldots,n-1$, where the coefficients
$\delta_{ij}$, $\alpha_i$ and $\beta_{i-1}$ are as defined in the Lemma~\ref{lemma:kalm.3}
The elements which are not defined are irrelevant to further considerations,
and are set to be zeros.
\smallskip

Consider a cut matrix in CDW normal form.  Let   
$I_l=\{i_l,\ldots,j_l\}$, $1\le l\le k$, be its $k$ blocks with more than one element,
involved in the  representation  described in Corollary~\ref{cor:1}. These
blocks  have the following properties:
$2\le|I_1|\le|I_2|\le \ldots \le |I_k|$,
$\cup_{l=1}^{l=k}I_l=\{i_1,i_1+1,\ldots,n\}$ and
$i_l=j_{l-1}+1$, for $2\le l\le k$.
Clearly, the corresponding  cut-weight matrix   contains only $k$ non-zero elements
$d_{i_l,j_l}=1$, $1\le l\le k$, exactly one  for each block. 
\smallskip

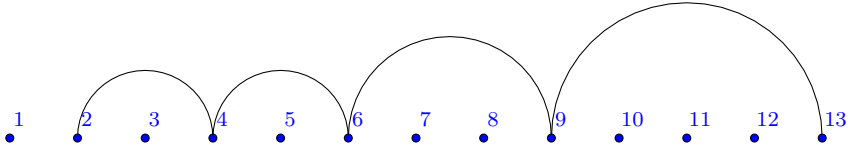
\begin{figure}[htb]
\definecolor{qqqqff}{rgb}{0.,0.,1.}
\begin{tikzpicture}[line cap=round,line join=round,>=triangle 45,x=0.9cm,y=0.9cm]
\clip(-3.0018181818181846,2.98) rectangle (10.5,6.66);
\draw [shift={(8.,4.)}] plot[domain=0.:3.141592653589793,variable=\t]({1.*2.*cos(\t r)+0.*2.*sin(\t r)},{0.*2.*cos(\t r)+1.*2.*sin(\t r)});
\draw [shift={(4.5,4.)}] plot[domain=0.:3.141592653589793,variable=\t]({1.*1.5*cos(\t r)+0.*1.5*sin(\t r)},{0.*1.5*cos(\t r)+1.*1.5*sin(\t r)});
\draw [shift={(2.,4.)}] plot[domain=0.:3.141592653589793,variable=\t]({1.*1.*cos(\t r)+0.*1.*sin(\t r)},{0.*1.*cos(\t r)+1.*1.*sin(\t r)});
\draw [shift={(0.,4.)}] plot[domain=0.:3.141592653589793,variable=\t]({1.*1.*cos(\t r)+0.*1.*sin(\t r)},{0.*1.*cos(\t r)+1.*1.*sin(\t r)});
\begin{scriptsize}
\draw [fill=qqqqff] (-2.,4.) circle (1.5pt);
\draw[color=qqqqff] (-1.8618181818181847,4.28) node {$1$};
\draw [fill=qqqqff] (-1.,4.) circle (1.5pt);
\draw[color=qqqqff] (-0.8618181818181847,4.28) node {$2$};
\draw [fill=qqqqff] (0.,4.) circle (1.5pt);
\draw[color=qqqqff] (0.1381818181818153,4.28) node {$3$};
\draw [fill=qqqqff] (1.,4.) circle (1.5pt);
\draw[color=qqqqff] (1.1381818181818153,4.28) node {$4$};
\draw [fill=qqqqff] (2.,4.) circle (1.5pt);
\draw[color=qqqqff] (2.1381818181818155,4.28) node {$5$};
\draw [fill=qqqqff] (3.,4.) circle (1.5pt);
\draw[color=qqqqff] (3.1381818181818155,4.28) node {$6$};
\draw [fill=qqqqff] (4.,4.) circle (1.5pt);
\draw[color=qqqqff] (4.138181818181815,4.28) node {$7$};
\draw [fill=qqqqff] (5.,4.) circle (1.5pt);
\draw[color=qqqqff] (5.138181818181815,4.28) node {$8$};
\draw [fill=qqqqff] (6.,4.) circle (1.5pt);
\draw[color=qqqqff] (6.138181818181815,4.28) node {$9$};
\draw [fill=qqqqff] (7.,4.) circle (1.5pt);
\draw[color=qqqqff] (7.198181818181816,4.28) node {$10$};
\draw [fill=qqqqff] (8.,4.) circle (1.5pt);
\draw[color=qqqqff] (8.198181818181814,4.28) node {$11$};
\draw [fill=qqqqff] (9.,4.) circle (1.5pt);
\draw[color=qqqqff] (9.198181818181816,4.28) node {$12$};
\draw [fill=qqqqff] (10.,4.) circle (1.5pt);
\draw[color=qqqqff] (10.198181818181816,4.28) node {$13$};
\end{scriptsize}
\end{tikzpicture}
\caption{Illustration for a block representation of a $12\times12$ cut matrix in CDW normal 
form with the blocks $\{2,3\}$, $\{4,5\}$, $\{6,7,8\}$, and $\{9,10,11,12\}$.}
\label{fig:matrix12}
\end{figure}

Next we will represent an   $n\times n$ cut matrix in CDW normal form by 
 a directed graph with $n+1$ nodes on a line, by means of the
cut-weight matrix, as follows. Let the nodes be labelled by $\{1,2,\ldots,n+1\}$,
increasing from the left to the right.
For each  non-zero entry $d_{i_1,j_1}=1$, $1\le l\le k$, of  the cut-weight matrix     we
introduce an edge that connects nodes $i_1$ and $j_1+1$ and is directed from $i_1$ to
$j_1+1$, hence from the left to the right; see Figure~\ref{fig:matrix12} for an illustration.
Let $k$ be the vertex with the smallest index haveing a positive degree. Then
the degree of $k$ equals $1$, i.e.  $deg(k)=1$ 
and every node $i\in \{k+1,\ldots,n\}$  has degree $0$ or
$2$,   whereas the degree of node $n+1$ equals $1$.
Furthermore notice that for every  directed edge   $(i,k)$  in this graph  $i+1<k$
holds. For such an edge we will say that it enters node $k$ and leaves node $i$. 
Finally notice that   if there is an  edge entering  a node $k\le n-1$,  then the 
edge leaving node $k$ is at least as long as the edge entering $k$, where
the length of an edge $(i,k)$ is given as $k-i$, for $i+1 < k$. 
Next we define a so-called multi-cut graph. 

\begin{definition}
A directed graph $G=(V,E)$ with node set
$V=\{1,2,\ldots,n+1\}$ and edge set $E=\{ (i_p,i_{p+1})\colon 1\le p\le |E|\}$,
where all edges are  directed from the node with the smaller index to the node with
larger index, with the follwowing properties. 
\begin{itemize}
\item[(a)] Let $i_1:=k$ be the vertex with the smallest index such that $deg(k)>0$. 
Then $deg(k)=1$, $deg(v)\in \{0,2\}$, for $v\in
  \{k+1,k+2,\ldots,n\}$ and  $deg(n+1)=1$ and $deg(v)\in \{0,2\}$, where $deg(v)$ denotes the degree of $v$ in $G$. 
\item[(b)] The length of every edge $(i_p,i_{p+1})\in E$ is not smaller than $2$,
  i.e. $i_{p+1}-i_p\ge 2$, for all  $1\le p \le |E|$. 
\item[(c)] For all vertices $i_{p+1}-i_p\le i_{p+2}-i_{p+i}$ holds for all
  $1\le p\le |E|-1$.
\end{itemize}
is called a {\sl cut-weight graph}.
\end{definition} 
 
It is straightforward to see that the symmetric matrix $D=(d_{ij})$ with
entries $d_{i_p,i_{p+1}}=1$ for $1\le p\le |E|$ and $d_{ij}=0$ otherwise,  for all
$i<j$,
is the cut-vertes matrix $D(C)=:D$ of an $n\times n$  block matrix $C$ in CDW normal
form. $C$ has $k-1+|E|$ blocks which are  given as 
$I_j=\{j\}$, for $j< k=i_1$, and  
$I_{k-1+t}=\{i_t,i_t+1,\ldots,i_{t+1}-1\}$, for $1\le t\le |E|$. 

So  to every  block matrix in CDW
normal form a cut-weight graph can be associated, and vice-cersa, to every cut-weight  graph a  block matrix in CDW
normal form can be associated,  as above.   
\smallskip

Consider now a conic combination $A=\sum_{p=1}^q\alpha_pA_p$ of (Kalmanson and Robinson) 
matrices $A_p$, for $1\le p\le q$. Clearly the cut-weight matrix of the linear (conic) combination can
be obtained as a  conic combination of the cut-weight  matrices of the
summands, i.e.\ $D(A)= \sum_{p=1}^q\alpha_p D(A_p)$. If $A_p$ are block matrices in
CDW normal form   and $\alpha_p\in \nz$, for all $1\le p\le q $, then the entries of $D(A)$
are natural numbers and $A$ can be represented as a directed multigraph with $n+1$ nodes on a
line  by means of its
cut-weight matrix $D(A)$. More precisely for every  two nodes $i$ and $k+1$, $1\le i\le
n-2$ and $i+2\le k+1\le n+1$,  the number of directed edges $(i,k+1)$  equals
 $\sum_{p=1}^q \alpha_p d^{(p)}_{i,k}$  where $d^{(p)}_{i,k}$ is
the corresponding entry of $D(A_p)$, for  $1\le p\le q$. Consequently for each
edge of length $x$ entering a node $k$, there is one edge of length at least $x$
leaving node $k$.  Let us denote by $E^-(k+1,x)$ and $E^+(k+1,x)$ the number of
edges of length at least $k$ entering or leaving the node $k+1$, respectively.
Then,  clearly  $E^-(k+1,x)\le E^+(k+1,x)$, holds  for $3\le k+1\le n-1$, $2\le
x\le k$, and by
considering that $E^-(k+1,x)= \sum_{i=1}^{k+1-x}d_{ik}$ and
$E^+(k+1,x)=\sum_{j=k+1+x-1}^n d_{k+1,j}$,  we get:
\begin{equation} \sum_{i=1}^{k+1-x}d_{ik}\le \sum_{j=k+1+x-1}^n d_{k+1,j}\, ,
  \mbox{ for $2\le k+1\le n-2$ and $2\le x\le k$.}  \label{upperbound:cutweights} \end{equation}

It turns out that these inequalities are not only necessary, but also
sufficient,  for the cut-weight matrix of a conic combination $A=\sum_{p=1}^q\alpha_pA_p$, 
 where $A_p$, $1\le p\le q$, are  block matrices in CDW normal form.

\begin{theorem}
A symmetric $n\times n$ Kalmanson matrix $C$ which is also a Robinson
matrix with  cut-weight matrix $D(C)=(d_{ij})$  can be represented as the sum
of a weak constant matrix and a conic
combination  of cut  matrices in CDW normal form 
if and only if the following inequalities hold
\begin{equation}
\sum_{i=1}^{l}d_{ik}\le\sum_{j=2k+1-l}^n d_{k+1,j}\label{eq:theorem}
\end{equation}
for  $k=2,\ldots,n-2$ and  $l=1,\ldots,k-1$.
\end{theorem}
The right-hand sum in \eqref{eq:theorem} is considered to be zero if $2k+1-l>n$. 
This is particular means that $d_{ik}=0$ for $k=\lceil n/2 \rceil,\ldots,n-1$ and $i=1,\ldots,2k-n$.
\smallskip

\proof 
Let $A=\sum_{p=1}^q \alpha_p A_p$ be an $n\times n$  matrix which is a conic
combination of $n\times n$ cut matrices
$A_p$ in CDW normal form,  $1\le p\le q$ 
Let $D(A)=(d_{ij})$ be  the 
 cut-weight matrix of $A$. 
Assume for simplicity (and without loss of generality) that the weight
coefficients $\alpha_p$, $1\le p\le q$, are natural numbers. Then, since
$D(A)=\sum_{p=1}^q D(A_p)$,   the entries $d_{lk}$, $1\le l,k\le n$, of the cut-weight matrix
$D(A)$ are also natural numbers.  
Inequality \eqref{eq:theorem} follows immediately from the inequality
\eqref{upperbound:cutweights} by setting $l:=k+1-x$.
\smallskip

Assume now that the entries of  the  cut-matrix $D(C)$ of  an integer Kalmanson
and symmetric matrix $C$ satisfy the inequalities \eqref{eq:theorem}. 
Similarly as in the discussion preceding the theorem we build an auxiliary
directed multigraph with $n+1$ nodes
$\{1,2,\ldots,n+1\}$  and  $d_{ij}$ edges directed from 
 $i$ to node $j+1$   for each $d_{ij}>0$ and $j> i$.  We refer to this
 multigraph as the \emph{cut-weight multigraph} of matrix $C$. 

We build a conic representation of $C$ as follows.
We start with node $n+1$ in the cut-weight multigraph, and build a  path from
the right to the left by choosing in the first step 
an (reversed) edge $(i_{k-1},i_k=n+1)$ with the largest length and then in each
following step $p> 1$  a
longest edge $(i_{k-p},i_{k-p+1})$ with length $i_{k-p+1}-i_{k-p}\le i_{k-p+2}-i_{k-p+1}$, as long
as  there is such an edge.   Let the  constructed path
be $P_1:=(i_1,i_2,\ldots, i_k=n+1)$ where  $1\le i_1<i_2<\ldots< i_k=n+1$. 
Consider now the graph $G=(V,E)$ with vertex set $V=\{1,2,\ldots,n+1\}$ and
edge set $E=\{(i_p,i_{p+1})\colon 1\le p\le k-1\}$. 
This is clearly a cut-weight graph and hence it can be associated to a cut
matrix  in CDW normal form as described above. We denote this matrix by
$A_1$. Let $\alpha_1$ be the number minimum multiplicity of the edges of the
path $P_1$. We remove $\alpha_1$ copies of $P_1$ from the cut-weight multi
graph and set $d_{i_p,i_{p+1}-1}:= d_{i_p,i_{p+1}-1}-\alpha_i$, for $1\le p\le
k-1$. We show that \eqref{eq:theorem} remains fulfilled after this update. 
Let us work with  the following equivalent formulation  of \eqref{eq:theorem}
\begin{equation} E^-(k+1,x)\le E^+(k+1,x)\,
  ,  \mbox{for  $k=2,\ldots,n-2$ and  $x=2,\ldots,k$,} \label{eq:theorem1}\end{equation}
which can be read as ``the number of edges of length at least $x$ entering
node $k+1$ does not exceed the number of edges of length at least $x$
leaving node $k+1$''. The update of the coefficients $d_{ij}$ can only affect
inequality \eqref{eq:theorem1} for indices $k$ such that $k+1$ is an endpoint of an
edge in $P_1$, i.e.\ $k+1\in \{i_1,i_2,\ldots, i_{k-1}\}$. For $k+1\in
\{i_2,\ldots, i_{k-1}\}$
the update of $d_{ij}$ results in  subtracting  $\alpha_i$ from
both sides of \eqref{eq:theorem1}, and hence it does not affect the validity of
the inequality. It remains the case  $k+1=i_1$ for $i_1>2$. 
Let $\bar{E}^-(i_1,x)$,  $\bar{E}^+(i_1,x)$ be the values of
$E^-(i_1,x)$ and $E^+(i_1,x)$, after the update of the coefficients $d_{ij}$,
respectively. Thus we have to
show that $\bar{E}^-(i_1,x)\le \bar{E}^+(i_1,x)$ holds. 
Let $l_1:=i_2-i_1$ be the length of the first edge in $P_1$. 
If $x>l_1$ than $\bar{E}^+(i_1,x)=E^+(i_1,x)$ and
$\bar{E}^-(i_1,x)=E^-(i_1,x)$, and since $E^-(i_1,x)\le E^+(i_1,x)$ holds, 
 there is nothing to show in this case. 
If $x\le l_1$, $\bar{E}^+(i_1,x)=E^+(i_1,x)-\alpha_1$ and
$\bar{E}^-(i_1,x)=E^-(i_1,x)$.
 Notice that $E^-(i_1,x)=E^-(i_1,x+1)$, otherwise  $P_1$ could have been
 prolonged beyong $i_1$. Moreover, $E^-(i_1,x+1)\le E^+(i_1,x+1)$ due to the
 assumption of the theorem,  and
 $E^+(i_1,x+1)\le E^+(i_1,x)-\alpha_1=\bar{E}^+(i_1,x)$ because there are at
 least $alpha_1$ edges of length $x$ leaving $i_1$. By putting things together
 we get the required inequality   $\bar{E}^-(i_1,x)\le \bar{E}^+(i_1,x)$.

The path construction and the corresponding  update of the coefficient $d_{ij}$
can be then inductively repeated as long as possible,  while \eqref{eq:theorem}
remains  an invariant during this process  and in every step $i$, $i\in \nz$,   a cut matrix
$A_i$ in CDW normal form is identified. $A_i$  corresponds to the path $P_i$ constructed  in
the $i$-th step. If $\alpha_i$ is the minimum multiplicity of
the edges in $P_i$ then $\alpha_iP_i$ is a summand of the required conic
combination. This process is finite because every step remove at least one edge
from the original cut-edge multigraph. The process terminates when there are no
edges entering node $n+1$ any more, say after $t$ steps.  We claim that after
the $t$-th step,              there are no more
edges in the cut-weight multigraph at all.  This means that    the actual coefficients $d_{ij}$
fulfill
     $d_{ij}=0$ for all $i<j$, which  implies that the matrix $C$ is 
     transformed into a weak contant matrix $Z$ by subtracting $\sum_{p=1}^t
     \alpha_pA_p$, and thus   $C=Z+\sum_{p=1}^t
     \alpha_pA_p$ holds for the original matrix $C$. 

Now let us prove the claim. Assume by contradiction that after $t$
transformation steps there is no edge entering node $n+1$ in the cut-weight
multigraph while there still at least one edge in it. Let $j$ be the largest
node index such that there is an edge entering $j$. Then $j\le n$ according to
our assumption. The inequalities  \eqref{eq:theorem1} have to be fulfilled
because they are an
invariant of the transformation process. In particular, $1\le E^-(j,1)\le
E^+(j,1)$ must hold. This implies the existece of an edge leaving $j$, and
hence entering some node with an index strictly larger than $j$. This
contradicts the choice of $j$ and completes the proof of the claim.     
\qed

\paragraph{An illustrative example.}
Consider  the  matrix $C$ which is in the class of Kalmanson and Robinson  matrices: 
\[
C=\left( \begin{array}{cccccc}
0&1&2&3&3&3\\
1&0&2&3&3&3\\
2&2&0&2&3&3\\
3&3&2&0&2&2\\
3&3&3&2&0&1\\
3&3&3&2&1&0\\
\end{array} \right)
\]
We first illustrate the proofs of Lemmas \ref{lemma:kalm.2} and \ref{lemma:kalm.3} and show how 
to represent $C$ as sum of  a conic combinatioon  of cut matrices $A^{kl}$ with
a weak
constant matrix.
 
Note that for matrix $C$ there is only one strict inequality in system \eqref{kalmx.c1}
$c_{25}+c_{34}<c_{24}+c_{35}$,
and three strict inequalities in system \eqref{kalmx.c2}:
$c_{i1}+c_{i+1,6}<c_{i6}+c_{i+1,1}$ with $i=2,3,4$.
 
We first eliminate the strict inequality $c_{25}+c_{34}<c_{24}+c_{35}$ by
subtracting from $C$ the cut matrix $A^{34}$ multiplied by
$\delta_{34}=c_{24}+c_{35}-c_{25}-c_{34}=d_{24}=1$. 
The  transformed matrix $C'=C-A^{34}$ is given as follows. 
\[
C'=C-A^{34}=\left( \begin{array}{cccccc}
0&0&1&2&2&2\\
0&0&1&2&2&2\\
1&1&0&2&2&2\\
2&2&2&0&1&1\\
2&2&2&1&0&0\\
2&2&2&1&0&0\\
\end{array} \right)
\]
Next we set  $\alpha_2=c'_{31}-c'_{21}=1$, $\beta_2=c'_{26}-c'_{36}=0$,
therefore in the next  
 transformation step the cut matrix $A^{12}$ is subtracted and we get:
\[
C-A^{34}-A^{12}=\left( \begin{array}{cccccc}
0&0&0&1&1&1\\
0&0&0&1&1&1\\
0&0&0&1&1&1\\
1&1&1&0&0&0\\
1&1&1&0&0&-1\\
1&1&1&0&-1&0\\
\end{array} \right)
\]
In the next step we set $\alpha_3=1$, $\beta_3=1$, and subtract $A^{13}+A^{46}$
from the current matrix $C$ to obtain 
\[
C-A^{34}-A^{12}-A^{13}-A^{46}=\left( \begin{array}{cccccc}
0&-1&-1&-1&-1&-1\\
-1&0&-1&-1&-1&-1\\
-1&-1&0&-1&-1&-1\\
-1&-1&-1&0&-1&-1\\
-1&-1&-1&-1&0&-2\\
-1&-1&-1&-1&-2&0\\
\end{array} \right)
\]
Finally we set $\alpha_4=c^{\prime}_{51}-c^{\prime}_{41}=0$ and
$\beta_4=c^{\prime}_{4,6}-c^{\prime}_{56}=1$ and subtract $A^{56}$ from the actual matrix $C$
to obtain a weak constant matrix 
\[
C-A^{34}-A^{12}-A^{13}-A^{46}-A^{56}=\left( \begin{array}{cccccc}
0&-2&-2&-2&-2&-2\\
-2&0&-2&-2&-2&-2\\
-2&-2&0&-2&-2&-2\\
-2&-2&-2&0&-2&-2\\
-2&-2&-2&-2&0&-2\\
-2&-2&-2&-2&-2&0\\
\end{array} \right)
\]
\smallskip

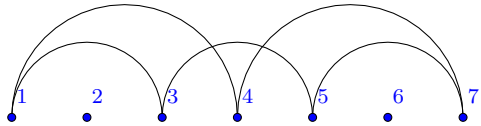
\begin{figure}[htb]
\begin{center}

\definecolor{qqqqff}{rgb}{0.,0.,1.}
\begin{tikzpicture}[line cap=round,line join=round,>=triangle 45,x=1.0cm,y=1.0cm]
\clip(-2.98,3.42) rectangle (4.5,5.92);
\draw [shift={(-1.,4.)}] plot[domain=0.:3.141592653589793,variable=\t]({1.*1.*cos(\t r)+0.*1.*sin(\t r)},{0.*1.*cos(\t r)+1.*1.*sin(\t r)});
\draw [shift={(-0.5,4.)}] plot[domain=0.:3.141592653589793,variable=\t]({1.*1.5*cos(\t r)+0.*1.5*sin(\t r)},{0.*1.5*cos(\t r)+1.*1.5*sin(\t r)});
\draw [shift={(1.,4.)}] plot[domain=0.:3.141592653589793,variable=\t]({1.*1.*cos(\t r)+0.*1.*sin(\t r)},{0.*1.*cos(\t r)+1.*1.*sin(\t r)});
\draw [shift={(3.,4.)}] plot[domain=0.:3.141592653589793,variable=\t]({1.*1.*cos(\t r)+0.*1.*sin(\t r)},{0.*1.*cos(\t r)+1.*1.*sin(\t r)});
\draw [shift={(2.5,4.)}] plot[domain=0.:3.141592653589793,variable=\t]({1.*1.5*cos(\t r)+0.*1.5*sin(\t r)},{0.*1.5*cos(\t r)+1.*1.5*sin(\t r)});
\begin{scriptsize}
\draw [fill=qqqqff] (-2.,4.) circle (1.5pt);
\draw[color=qqqqff] (-1.86,4.28) node {$1$};
\draw [fill=qqqqff] (-1.,4.) circle (1.5pt);
\draw[color=qqqqff] (-0.86,4.28) node {$2$};
\draw [fill=qqqqff] (0.,4.) circle (1.5pt);
\draw[color=qqqqff] (0.14,4.28) node {$3$};
\draw [fill=qqqqff] (1.,4.) circle (1.5pt);
\draw[color=qqqqff] (1.14,4.28) node {$4$};
\draw [fill=qqqqff] (2.,4.) circle (1.5pt);
\draw[color=qqqqff] (2.14,4.28) node {$5$};
\draw [fill=qqqqff] (3.,4.) circle (1.5pt);
\draw[color=qqqqff] (3.14,4.28) node {$6$};
\draw [fill=qqqqff] (4.,4.) circle (1.5pt);
\draw[color=qqqqff] (4.14,4.28) node {$7$};
\end{scriptsize}
\end{tikzpicture}
\end{center}
\caption{ Cut weight graph for the illustrative example}
\label{fig:illustration}
\end{figure}

The cut-weight matrix $D(C)$ contains five non-zero entries corresponding to
the coefficients $\delta_{34}$, $\alpha_2$, $\alpha_3$, $\beta_3$ and $\beta_4$ above:
\[
D(C) =\left( \begin{array}{cccccc}
0&1&1&0&0&0\\
0&0&0&0&0&0\\
0&0&0&1&0&0\\
0&0&0&0&0&1\\
0&0&0&0&0&1\\
0&0&0&0&0&0\\
\end{array} \right)
\]
The corresponding cut-weight multigraph is depicted in   Figure~\ref{fig:illustration}.
It can be easily seen that the entries of $D(C)$ fulfill the
inequalities~\ref{eq:theorem} 
and hence $C$ can be represented as the sum of weak constant
matrix with a conic combination of block matrices in CDW normal form. 
The block matrices $A_1$ and $A_2$ in CDW normal form correspond to the paths $P_1=(1,4,7)$ and
$P_2=(1,3,5,7)$  and hence, $A_1$  has two blocks $\{1,2,3\}$ and $\{4,5,6\}$,
and  $A_2$ 
has three blocks $\{1,2\}$, $\{3,4\}$, and $\{5,6\}$.

\section{Conclusions}\label{conclu:sec}
In this paper we introduced two new polynomially solvable special cases of the
QAP. 
We call the first one {\sl the down-benevolent QAP}; this is a  $QAP(A,B)$ where $A$ 
is both a Kalmanson and a Robinson matrix and
$B$ is  a down-benevolent Toeplitz matrix, and it is solved to otimality by the
identity   permutation. This new special case is a generaliation of two other special cases of the QAP
 known in the literature: (a) the $QAP(A,B)$ with $A$ being a Kalmanson
matrix and $B$ being a DW Toeplitz matrix~\cite{DeWo1998}, and (b) the $QAP(A,B)$ with $A$
being a Robinson matrix and $B$ being a simple Toeplitz
matrix~\cite{Laurent2015}. 

We call the second new special case {\sl the up-benevolent QAP}; this is a QAP(A,B) where $A$ is a
PS monotone Anti-Monge matrix and $B$ is an up-benevolent Toeplitz matrix,  and
it is solved to optimality by the identity permutation. This new special case
is a generalization of 
another special case of the QAP  known in the literature, namely  the
$QAP(A,B)$ where $A$ is  a symmetric monotone Anti-Monge matrix   and $B$ is
an up-benevolent Toeplitz matrix~\cite{BCRW1998}.

Further we  introduce a new class of specially structured matrices. 
A matrix belongs to this  class if it can be represented as the sum of a
weakly constant matrix and a conic combination of cut matrices in CDW normal
form.   The matrices of this class build a strict subclass of matrices which are both Robinson and Kalmanson matrices. 
It follows from a result in \cite{Cela2015} that the  $QAP(A,B)$ is solved to
optimality by the identity permutation if  $A$ belongs
to the newly introduced class of matrices and $B$ is a symmetric monotone anti-Monge matrix 
 
The new class of matrices and the down-benevolent QAP
 lead to another new  polynomially solvable special case of the
QAP, namely a combined special case:
 the  $QAP(A,B)$ where   $A$ is a conic combination  of cut matrices in CDW normal
form and $B$  is  a conic combination of a monotone  anti-Monge 
matrix 
and a down-benevolent Toeplitz matrix, is solved to optimality by the identity
permutation. 

The combined special case mentioned above gives rise to an interesting and
non-trivial 
question related to  the recognition of conic combination of cut matrices in CDW normal
form: Given an $n\times n$ matrix
$A$, $n\in \nz$, decide whether  $A$  can be represented as the sum of a
weak constant matrix and a conic combination of cut matrices in CDW normal
form. 
We show that this decision problem can be solved efficiently by computing $O(n^2)$
so-called cut-weights and checking whether these weights fulfill $O(n^2)$
linear inequalities.  

Notice that both the   monotone anti-Monge matrices and the down-benevolent Toeplitz
matrices are defined in terms of linear inequalities. Therefore
 simple linear programming techniques can be used to recognize whether a given
 symmetric matrix $B$  can be 
represented/approximated as a conic combination of two  matrices $B_1$ and $B_2$
where $B_1$ is a monotone anti-Monge matrix  and $B_2$ is  a down-benevolent Toeplitz
matrices.  Thus for a given instance of the  $QAP(A,B)$ it can be efficiently
checked whether it is an instance of the new combined  special case
introduced in this paper. 

A more general and challenging question to be
considered for future research   is the recognition of the
so-called permuted combined special case. For a given an instance $QAP(A,B)$ decide
whether a)  there exists a permutation $\phi$ of the rows and the columns of $A$
such that the matrix resulting after permuting $A$ according to $\phi$ is the
sum of a weak constant matrix and a
conic combination of  cut matrices in CDW normal
form, and b) there exists a permutation $\psi$ of the rows and the columns of $B$
such that the matrix resulting after permuting $B$ according to $\psi$ can be
represented as the sum of a monotone anti-Monge matrix and a down-benevolent Toeplitz
matrix.

 Moreover it would be interesting to investigate whether  the  combined special case of the QAP can be used to
compute good lower bounds and/or heuristic solutions  for the general
problem. The idea is to  ``approximate'' the coefficient matrices $A$ and
$B$ of a given  instance $QAP(A,B)$ by some  matrices $A'$ and $B'$,
respectively,  such that the 
$QAP(A',B')$ is an instance of the combined special case. Then, if $A'$  and
$B'$ are chosen ``appropriately'', the optimal
solution of $QAP(A',B')$  and its optimal value could serve as a heuristic
solution  and/or a lower bound for the $QAP(A,B)$, respectively. Clearly, the crucial 
 part is to find out what ``approximate'' and ``appropriately''
should mean.  This is definitely a challenging issue but it could well 
lead to a new direction of research on the QAP.    
 
\paragraph{Acknowledgements.}
Vladimir Deineko acknowledges support
by Warwick University's Centre for Discrete Mathematics and Its Applications (DIMAP).
Gerhard Woeginger acknowledges support
by the Zwaartekracht NETWORKS grant of NWO,
and by the Alexander von Humboldt Foundation, Bonn, Germany.

\newpage

\begin{appendix}
\section{Instances used in illustrations}
\begin{figure}[ht]
\unitlength=1cm
\begin{center}
\begin{picture}(10.5,6)
\put(7,0.5)
{\begin{picture}(3.8,4)
\includegraphics[scale=1]{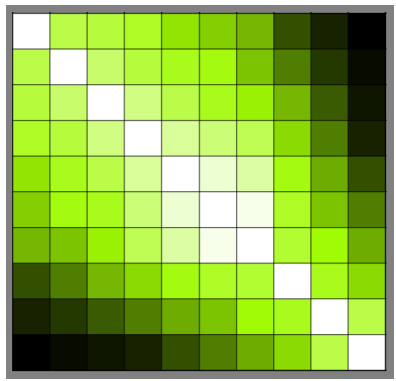}
\end{picture} }
\put(-2.5,2.9){
{
\begin{tabular}{|c|c|c|c|c|c|c|c|c|c|}
\hline 
0 & 16 & 17 & 19 & 24 & 26 & 28 & 37 & 41 & 44\\
\hline
16 & 0 & 13 & 17 & 20 & 21 & 27 & 33 & 39 & 43\\
\hline
17 & 13 & 0 & 11 & 16 & 20 & 23 & 28 & 36 & 42\\
\hline
19 & 17 & 11 & 0 & 9 & 12 & 15 & 25 & 33 & 41\\
\hline
24 & 20 & 16 & 9 & 0 & 4 & 8 & 21 & 29 & 37\\
\hline
26 & 21 & 20 & 12 & 4 & 0 & 2 & 19 & 27 & 33\\
\hline
28 & 27 & 23 & 15 & 8 & 2 & 0 & 18 & 22 & 29\\
\hline
37 & 33 & 28 & 25 & 21 & 19 & 18 & 0 & 20 & 25\\
\hline
41 & 39 & 36 & 33 & 29 & 27 & 22 & 20 & 0 & 16\\
\hline
44 & 43 & 42 & 41 & 37 & 33 & 29 & 25 & 16 & 0\\
\hline
\end{tabular}
} 
 }
\end{picture}
\end{center}
\vspace{-1cm}
\caption{Robinsonian dissimilarity visualised in Fig.~\ref{fig:RobToeplitz} A}
\label{fig:Instance1.1}
\end{figure}

\begin{figure}[ht]
\unitlength=1cm
\begin{center}
\begin{picture}(10.5,5.8)
\put(7,0.5)
{\begin{picture}(3.8,4)
\includegraphics[scale=1]{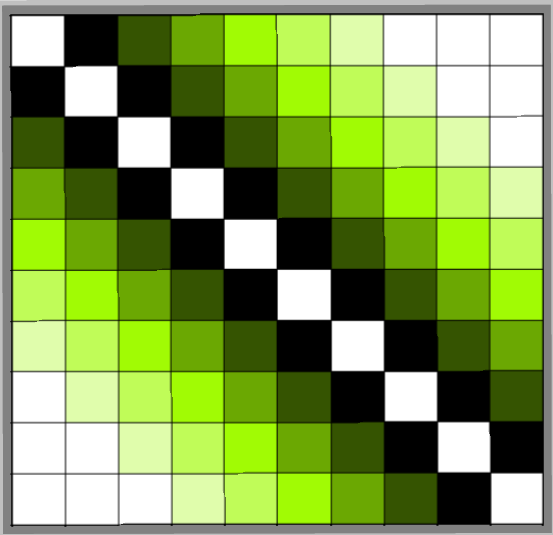}
\end{picture} }
\put(-2.5,2.9){
{
\begin{tabular}{|c|c|c|c|c|c|c|c|c|c|}
\hline 
0 & 55 & 45 & 35 & 25 & 15 & 5 & 5 & 0 & 0\\
\hline
55 & 0 & 55 & 45 & 5 & 3 & 15 & 5 & 5 & 0\\
\hline
45 & 55 & 0 & 55 & 45 & 5 & 3 & 15 & 5 & 5\\
\hline
35 & 45 & 55 & 0 & 55 & 45 & 5 & 3 & 15 & 5\\
\hline
25 & 5 & 45 & 55 & 0 & 55 & 45 & 5 & 3 & 15\\
\hline
15 & 3 & 5 & 45 & 55 & 0 & 55 & 45 & 5 & 3\\
\hline
5 & 15 & 3 & 5 & 45 & 55 & 0 & 55 & 45 & 5\\
\hline
5 & 5 & 15 & 3 & 5 & 45 & 55 & 0 & 55 & 45\\
\hline
0 & 5 & 5 & 15 & 3 & 5 & 45 & 55 & 0 & 55\\
\hline
0 & 0 & 5 & 5 & 15 & 3 & 5 & 45 & 55 & 0\\
\hline
\end{tabular}
} 
 }
\end{picture}
\end{center}
\vspace{-1cm}
\caption{Simple Toeplitz matrix visualised in Fig.~\ref{fig:RobToeplitz} B}
\label{fig:Instance1.2}
\end{figure}

\begin{figure}[ht]
\unitlength=1cm
\begin{center}
\begin{picture}(10.5,5.8)
\put(7,0.5)
{\begin{picture}(3.8,4)
\includegraphics[scale=1]{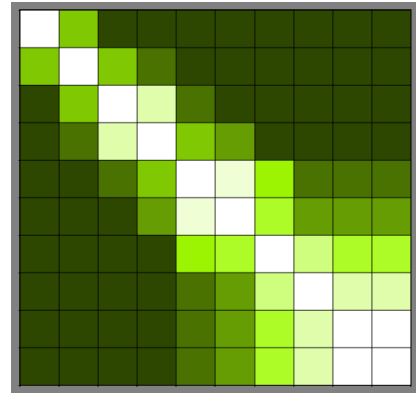}
\end{picture} }
\put(-2.5,2.9){
{
\begin{tabular}{|c|c|c|c|c|c|c|c|c|c|}
\hline 
0 & 7 & 10 & 10 & 10 & 10 & 10 & 10 & 10 & 10\\
\hline
7 & 0 & 7 & 9 & 10 & 10 & 10 & 10 & 10 & 10\\
\hline
10 & 7 & 0 & 2 & 9 & 10 & 10 & 10 & 10 & 10\\
\hline
10 & 9 & 2 & 0 & 7 & 8 & 10 & 10 & 10 & 10\\
\hline
10 & 10 & 9 & 7 & 0 & 1 & 6 & 9 & 9 & 9\\
\hline
10 & 10 & 10 & 8 & 1 & 0 & 5 & 8 & 8 & 8\\
\hline
10 & 10 & 10 & 10 & 6 & 5 & 0 & 3 & 5 & 5\\
\hline
10 & 10 & 10 & 10 & 9 & 8 & 3 & 0 & 2 & 2\\
\hline
10 & 10 & 10 & 10 & 9 & 8 & 5 & 2 & 0 & 0\\
\hline
10 & 10 & 10 & 10 & 9 & 8 & 5 & 2 & 0 & 0\\
\hline
\end{tabular}
} 
 }
\end{picture}
\end{center}
\vspace{-1cm}
\caption{Conic combination of block matrices in CDW normal form  visualised in Fig.~\ref{fig:BlockAntiM} A}
\label{fig:Instance2.1}
\end{figure}

\begin{figure}[ht]
\unitlength=1cm
\begin{center}
\begin{picture}(10.5,5.8)
\put(7.2,0.5)
{\begin{picture}(3.8,4)
\includegraphics[scale=0.95]{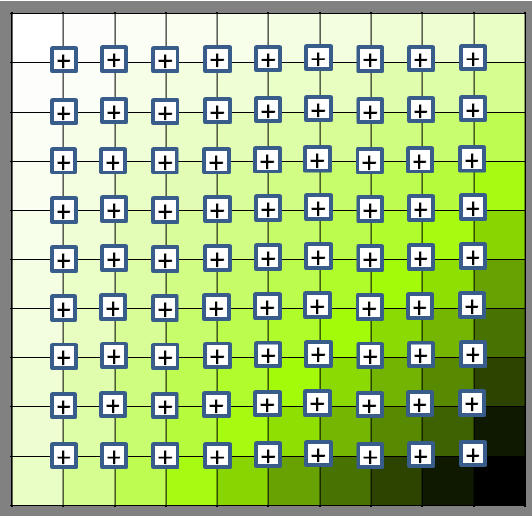}
\end{picture} }
\put(-2.7,2.9){
{
\begin{tabular}{|c|c|c|c|c|c|c|c|c|c|}
\hline 
1 & 3 & 6 & 7 & 10 & 13 & 14 & 18 & 19 & 24\\
\hline
3 & 5 & 11 & 14 & 18 & 23 & 26 & 33 & 36 & 45\\
\hline
6 & 11 & 17 & 24 & 30 & 36 & 42 & 51 & 57 & 69\\
\hline
7 & 14 & 24 & 31 & 40 & 49 & 59 & 69 & 77 & 93\\
\hline
10 & 18 & 30 & 40 & 49 & 60 & 71 & 83 & 94 & 114\\
\hline
13 & 23 & 36 & 49 & 60 & 71 & 84 & 97 & 110 & 134\\
\hline
14 & 26 & 42 & 59 & 71 & 84 & 97 & 112 & 126 & 152\\
\hline
18 & 33 & 51 & 69 & 83 & 97 & 112 & 127 & 143 & 170\\
\hline
19 & 36 & 57 & 77 & 94 & 110 & 126 & 143 & 159 & 187\\
\hline
24 & 45 & 69 & 93 & 114 & 134 & 152 & 170 & 187 & 215\\
\hline
\end{tabular}
} 
 }
\end{picture}
\end{center}
\vspace{-1cm}
\caption{Anti-Monge monotone matrix  visualised in Fig.~\ref{fig:BlockAntiM} B}
\label{fig:Instance2.2}
\end{figure}

\begin{figure}[ht]
\unitlength=1cm
\begin{center}
\begin{picture}(10.5,5.8)
\put(7.2,0.5)
{\begin{picture}(3.8,4)
\includegraphics[scale=0.95]{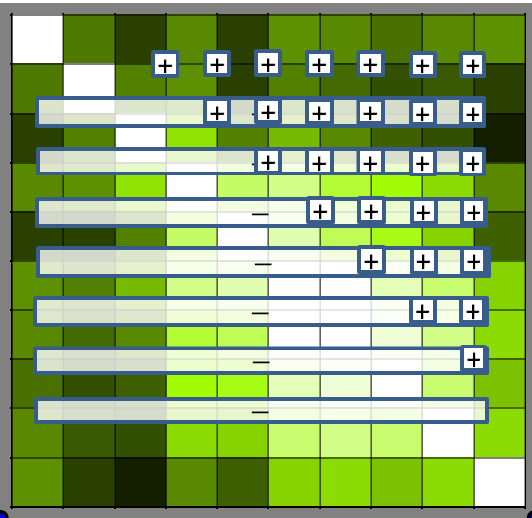}
\end{picture} }
\put(-2.7,2.9){
{
\begin{tabular}{|c|c|c|c|c|c|c|c|c|c|}
\hline 
0 & 47 & 54 & 45 & 54 & 44 & 45 & 48 & 45 & 44\\				
\hline				
47 & 0 & 46 & 44 & 54 & 46 & 49 & 52 & 51 & 54\\				
\hline				
54 & 46 & 0 & 34 & 46 & 39 & 45 & 50 & 52 & 58\\				
\hline				
45 & 44 & 34 & 0 & 19 & 15 & 25 & 31 & 35 & 45\\				
\hline				
54 & 54 & 46 & 19 & 0 & 10 & 21 & 29 & 36 & 50\\				
\hline				
44 & 46 & 39 & 15 & 10 & 0 & 0 & 9 & 18 & 36\\				
\hline				
45 & 49 & 45 & 25 & 21 & 0 & 0 & 5 & 15 & 35\\				
\hline				
48 & 52 & 50 & 31 & 29 & 9 & 5 & 0 & 18 & 38\\				
\hline				
45 & 51 & 52 & 35 & 36 & 18 & 15 & 18 & 0 & 35\\				
\hline				
44 & 54 & 58 & 45 & 50 & 36 & 35 & 38 & 35 & 0\\				
\hline
\end{tabular}
} 
 }
\end{picture}
\end{center}
\vspace{-1cm}
\caption{Kalmanson matrix  visualised in Fig.~\ref{fig:DeWo1998} A}
\label{fig:Instance3.1}
\end{figure}

\begin{figure}[ht]
\unitlength=1cm
\begin{center}
\begin{picture}(10.5,5.8)
\put(7.2,0.5)
{\begin{picture}(3.8,4)
\includegraphics[scale=0.95]{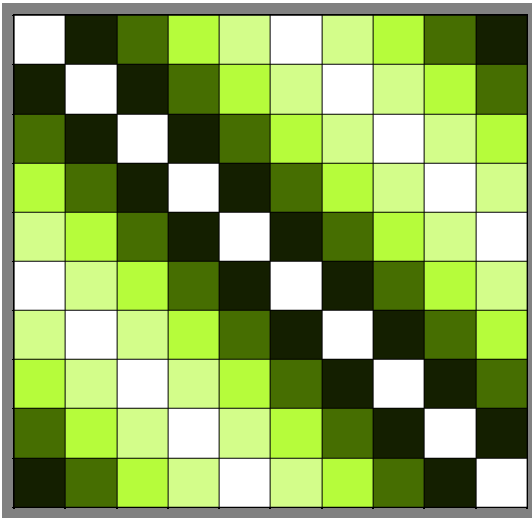}
\end{picture} }
\put(-2.7,2.9){
{
\begin{tabular}{|c|c|c|c|c|c|c|c|c|c|}
\hline 
0 & 12 & 10 & 5 & 3 & 0 & 3 & 5 & 10 & 12\\
\hline
12 & 0 & 12 & 10 & 5 & 3 & 0 & 3 & 5 & 10\\
\hline
10 & 12 & 0 & 12 & 10 & 5 & 3 & 0 & 3 & 5\\
\hline
5 & 10 & 12 & 0 & 12 & 10 & 5 & 3 & 0 & 3\\
\hline
3 & 5 & 10 & 12 & 0 & 12 & 10 & 5 & 3 & 0\\
\hline
0 & 3 & 5 & 10 & 12 & 0 & 12 & 10 & 5 & 3\\
\hline
3 & 0 & 3 & 5 & 10 & 12 & 0 & 12 & 10 & 5\\
\hline
5 & 3 & 0 & 3 & 5 & 10 & 12 & 0 & 12 & 10\\
\hline
10 & 5 & 3 & 0 & 3 & 5 & 10 & 12 & 0 & 12\\
\hline
12 & 10 & 5 & 3 & 0 & 3 & 5 & 10 & 12 & 0\\				
\hline
\end{tabular}
} 
 }
\end{picture}
\end{center}
\vspace{-1cm}
\caption{Toeplitz matrix  visualised in Fig.~\ref{fig:DeWo1998} B}
\label{fig:Instance3.2}
\end{figure}

\begin{figure}[!ht]
\unitlength=1cm
\begin{center}
\begin{picture}(10.5,5.8)
\put(7.2,0.5)
{\begin{picture}(3.8,4)
\includegraphics[scale=0.95]{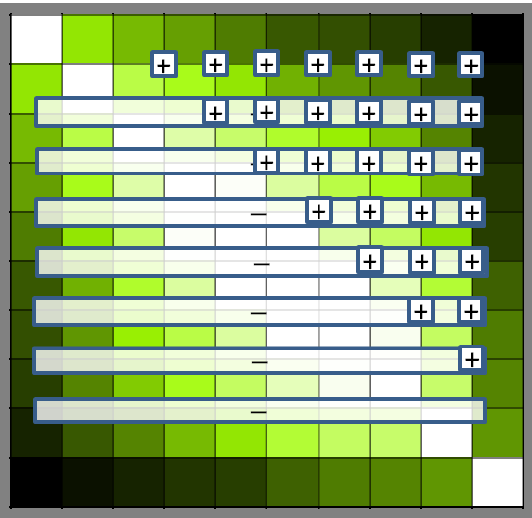}
\end{picture} }
\put(-2.7,2.9){
{
\begin{tabular}{|c|c|c|c|c|c|c|c|c|c|}
\hline 
0 & 31 & 36 & 39 & 43 & 47 & 48 & 50 & 53 & 57\\
\hline
31 & 0 & 21 & 26 & 31 & 37 & 40 & 42 & 47 & 55\\
\hline
36 & 21 & 0 & 10 & 17 & 24 & 30 & 34 & 42 & 53\\
\hline
39 & 26 & 10 & 0 & 1 & 11 & 21 & 26 & 36 & 51\\
\hline
43 & 31 & 17 & 1 & 0 & 0 & 11 & 18 & 31 & 50\\
\hline
47 & 37 & 24 & 11 & 0 & 0 & 0 & 8 & 23 & 46\\
\hline
48 & 40 & 30 & 21 & 11 & 0 & 0 & 2 & 18 & 43\\
\hline
50 & 42 & 34 & 26 & 18 & 8 & 2 & 0 & 17 & 42\\
\hline
53 & 47 & 42 & 36 & 31 & 23 & 18 & 17 & 0 & 40\\
\hline
57 & 55 & 53 & 51 & 50 & 46 & 43 & 42 & 40 & 0\\			
\hline
\end{tabular}
} 
 }
\end{picture}
\end{center}
\vspace{-1cm}
\caption{Kalmanson and Robinsonian matrix visualised in Fig.~\ref{fig:Malevolent} A}
\label{fig:Instance4.1}
\end{figure}

\begin{figure}[ht]
\unitlength=1cm
\begin{center}
\begin{picture}(10.5,5.8)
\put(7.2,0.5)
{\begin{picture}(3.8,4)
\includegraphics[scale=0.95]{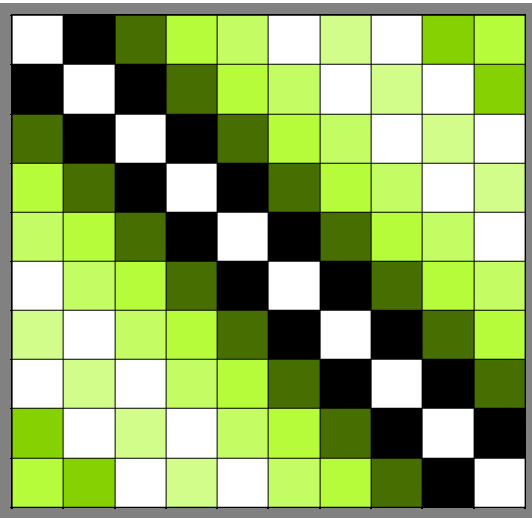}
\end{picture} }
\put(-2.7,2.9){
{
\begin{tabular}{|c|c|c|c|c|c|c|c|c|c|}
\hline 
0 & 28 & 20 & 10 & 8 & 0 & 6 & 0 & 15 & 10\\
\hline
28 & 0 & 28 & 20 & 10 & 8 & 0 & 6 & 0 & 15\\
\hline
20 & 28 & 0 & 28 & 20 & 10 & 8 & 0 & 6 & 0\\
\hline
10 & 20 & 28 & 0 & 28 & 20 & 10 & 8 & 0 & 6\\
\hline
8 & 10 & 20 & 28 & 0 & 28 & 20 & 10 & 8 & 0\\
\hline
0 & 8 & 10 & 20 & 28 & 0 & 28 & 20 & 10 & 8\\
\hline
6 & 0 & 8 & 10 & 20 & 28 & 0 & 28 & 20 & 10\\
\hline
0 & 6 & 0 & 8 & 10 & 20 & 28 & 0 & 28 & 20\\
\hline
15 & 0 & 6 & 0 & 8 & 10 & 20 & 28 & 0 & 28\\
\hline
10 & 15 & 0 & 6 & 0 & 8 & 10 & 20 & 28 & 0\\			
\hline
\end{tabular}
} 
 }
\end{picture}
\end{center}
\vspace{-1cm}
\caption{Malevolent matrix visualised in Fig.~\ref{fig:Malevolent} B}
\label{fig:Instance4.2}
\end{figure}

\begin{figure}[ht]
\unitlength=1cm
\begin{center}
\begin{picture}(10.5,5.8)
\put(7.2,0.5)
{\begin{picture}(3.8,4)
\includegraphics[scale=0.95]{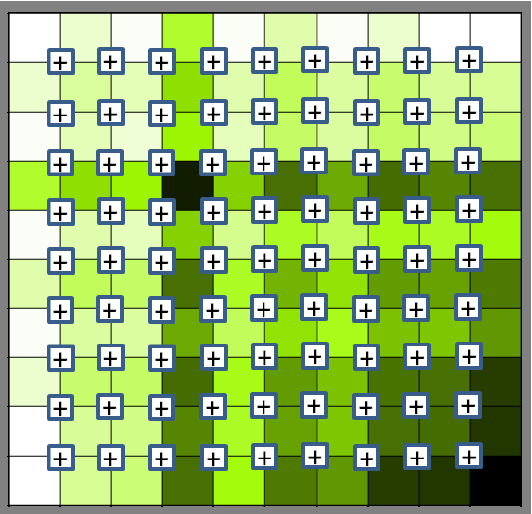}
\end{picture} }
\put(-2.7,2.9){
{
\begin{tabular}{|c|c|c|c|c|c|c|c|c|c|}
\hline 
1 & 15 & 4 & 60 & 4 & 25 & 4 & 15 & 1 & 1\\
\hline
15 & 29 & 21 & 79 & 24 & 47 & 28 & 42 & 30 & 31\\
\hline
4 & 21 & 13 & 73 & 20 & 44 & 28 & 44 & 35 & 39\\
\hline
60 & 79 & 73 & 133 & 83 & 110 & 94 & 111 & 104 & 110\\
\hline
4 & 24 & 20 & 83 & 33 & 62 & 47 & 66 & 62 & 69\\
\hline
25 & 47 & 44 & 110 & 62 & 91 & 78 & 98 & 96 & 107\\
\hline
4 & 28 & 28 & 94 & 47 & 78 & 65 & 87 & 86 & 99\\
\hline
15 & 42 & 44 & 111 & 66 & 98 & 87 & 109 & 110 & 124\\
\hline
1 & 30 & 35 & 104 & 62 & 96 & 86 & 110 & 111 & 126\\
\hline
1 & 31 & 39 & 110 & 69 & 107 & 99 & 124 & 126 & 141\\		
\hline
\end{tabular}
} 
 }
\end{picture}
\end{center}
\vspace{-1cm}
\caption{Anti-Monge matrix visualised in Fig.~\ref{fig:Bene} A}
\label{fig:Instance5.1}
\end{figure}

\begin{figure}[!ht]
\unitlength=1cm
\begin{center}
\begin{picture}(10.5,5.8)
\put(7.2,0.5)
{\begin{picture}(3.8,4)
\includegraphics[scale=0.95]{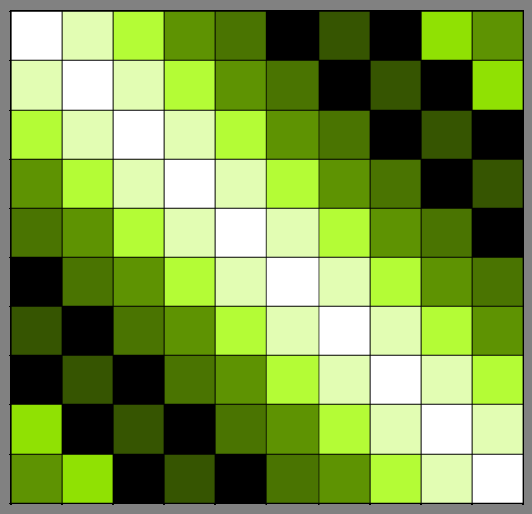}
\end{picture} }
\put(-2.7,2.9){
{
\begin{tabular}{|c|c|c|c|c|c|c|c|c|c|}
\hline 
0 & 5 & 13 & 23 & 25 & 33 & 27 & 33 & 18 & 23\\
\hline
5 & 0 & 5 & 13 & 23 & 25 & 33 & 27 & 33 & 18\\
\hline
13 & 5 & 0 & 5 & 13 & 23 & 25 & 33 & 27 & 33\\
\hline
23 & 13 & 5 & 0 & 5 & 13 & 23 & 25 & 33 & 27\\
\hline
25 & 23 & 13 & 5 & 0 & 5 & 13 & 23 & 25 & 33\\
\hline
33 & 25 & 23 & 13 & 5 & 0 & 5 & 13 & 23 & 25\\
\hline
27 & 33 & 25 & 23 & 13 & 5 & 0 & 5 & 13 & 23\\
\hline
33 & 27 & 33 & 25 & 23 & 13 & 5 & 0 & 5 & 13\\
\hline
18 & 33 & 27 & 33 & 25 & 23 & 13 & 5 & 0 & 5\\
\hline
23 & 18 & 33 & 27 & 33 & 25 & 23 & 13 & 5 & 0\\	
\hline
\end{tabular}
} 
 }
\end{picture}
\end{center}
\vspace{-1cm}
\caption{Benevolent Toeplitz matrix visualised in Fig.~\ref{fig:Bene} B}
\label{fig:Instance5.2}
\end{figure}

\end{appendix}

\end{document}